\def\versione{\label{vers.16}}

\documentclass[12pt]{article}

\usepackage{amssymb}
\usepackage{mathrsfs}
\usepackage{theorem}

\newtheorem{theorem}{Theorem}
\newtheorem{lemma}{Lemma}
\newtheorem{corollary}{Corollary}
\newenvironment{proof} {{\bf Proof.}}{\hfill \fbox{}\\ \smallskip}

\newcommand{\paragrafo}[1]{\section{#1}\setcounter{equation}{0}}

\makeatletter
\renewcommand{\@makefnmark}{\hbox
{\@textsuperscript{\normalfont(\@thefnmark)}}}
\makeatother
      
{\theorembodyfont{\rmfamily} \newtheorem{remark}{Remark}}

\newcommand{\reff}[1]{{\rm(\ref{#1})}}

\def\display{\displaylines}
\def\grande{\displaystyle}

\newcommand{\beginrighe}[1]{\begin{eqnarray}\label{#1}\begin{array}{c}\grande}
\def\endrighe{\end{array}\end{eqnarray}}

\def\C{\mathbb{C}}
\def\R{\mathbb{R}}

\def\de{\partial}

\def\essinf{\mathop{\rm ess\, inf}}
\def\esssup{\mathop{\rm ess\, sup}}

\def\a{\alpha}
\def\b{\beta}
\def\g{\gamma}

\def\d{\delta}
\def\D{\Delta}
\def\e{\varepsilon}
\def\th{\vartheta}

\def\h{\eta}
\def\l{\lambda}
\def\L{\Lambda}
\def\m{\mu}
\def\n{\nu}

\def\f{\varphi}

\def\O{\Omega}
\def\o{\omega}
\def\p{\psi}

\def\Cspt{\mathop{C}_{0}}

\def\tr{\mathop{\rm tr}\nolimits}
\def\arccot{\mathop{\rm arccot}\nolimits}

\def\Re{\mathop{\mathscr Re}\nolimits}
\def\Im{\mathop{\mathscr Im}\nolimits}
\def\elle{\mathop{\mathscr L}\nolimits}
\def\A{\mathop{\mathscr A}\nolimits}

\def\spt{\mathop{\rm spt}\nolimits}

\def\dive{\mathop{\rm div}\nolimits}

\def\leq{\leqslant}
\def\geq{\geqslant}
\def\lan{\langle}
\def\ran{\rangle}

\title{
Criteria for the $L^{p}$-dissipativity of systems
of second order differential equations
}
\author{A. Cialdea 
\thanks{Dipartimento di Matematica, Universit\`a della Basilicata, Viale 
dell'Ateneo Lucano 10, 85100, Potenza, Italy. \textit{email:}
cialdea@email.it.}  
\and V. Maz'ya
\thanks{Department of Mathematics, Ohio State University,
231 W 18th Avenue, Columbus, OH 43210, USA.
Department of Mathematical Sciences, M{\&}O Building,
University of Liverpool, Liverpool L69 3BX, UK.
\textit{email:} vlmaz@mai.liu.se.}
}
\date{}    

\begin{document}

    \maketitle

\bigskip

\versione
    {\small  {\bf Abstract.} We give complete algebraic
   characterizations of the 
    $L^{p}$-dissipativity of the Dirichlet problem for some
    systems of
    partial differential operators of the form 
    $\de_{h}(\A^{hk}(x)\de_{k})$, were $\A^{hk}(x)$ are $m\times m$
     matrices. 
     First, we determine the sharp angle of dissipativity
     for a general scalar operator with complex coefficients. 
    Next we prove that the two-dimensional elasticity 
    operator is $L^{p}$-dissipative if and only if 
    $$
    \left({1\over 2}-{1\over p}\right)^{2} \leq {2(\n-1)(2\n-1)\over
       (3-4\n)^{2}},
   $$
   $\n$ being the Poisson ratio. Finally we find a 
    necessary and sufficient algebraic condition for the $L^{p}$-dissipativity
    of the operator $\de_{h} (\A^{h}(x)\de_{h})$, where $\A^{h}(x)$ are 
    $m\times m$ matrices with complex $L^{1}_{\rm loc}$
    entries, and we describe the maximum angle of $L^{p}$-dissipativity
    for this operator.
    }

    \bigskip \bigskip

\paragrafo{Introduction}

Let $\O$ be a domain of $\R^{n}$ and
let $A$ be the operator
\begin{equation}
    A=\de_{h} (\A^{hk}(x)\de_{k})
    \label{eq:A}
\end{equation}
where  $\de_{k} =\de / \de x_{k}$ and $\A^{hk}(x)=\{a^{hk}_{ij}(x)\}$ are 
$m\times m$ matrices 
whose elements are complex locally integrable  functions defined 
in $\O$ 
$(1\leq i,j\leq 
m,\ 1\leq h,k \leq 2)$. Here and in the sequel we  adopt the 
summation convention and we put $p\in(1,\infty)$, $p'=p/(p-1)$.
By $\Cspt^{1}(\O)$ we denote the space of all the $C^{1}$
functions having 
compact support in $\O$.

Let $\elle$ be the sesquilinear form related to the operator $A$
$$
\elle(u,v) = \int_{\O}\lan \A^{hk}(x)\de_{k} u, \de_{h} v\ran\, dx.
$$
($\langle\cdot,\cdot \rangle$ denotes the scalar product in 
$\C^{m}$) defined in $(\Cspt^{1}(\O))^{m}\times (\Cspt^{1}(\O))^{m}$.
We consider $A$ as an operator acting from $(\Cspt^{1}(\O))^{m}$
to $((\Cspt^{1}(\O))^{*})^{m}$ through the relation
$$
\elle(u,v)=-\int_{\O}\lan Au, v\ran\, dx
$$
for any $u,v\in (\Cspt^{1}(\O))^{m}$. Here the integration is
understood in the sense of distributions.

Following \cite{cialmaz}, 
we say that the form $\elle$ is $L^{p}$-dissipative 
if
\begin{eqnarray}
    \Re {\mathscr L}(u, |u|^{p-2}u) \geq 0
    \qquad \hbox{\rm if}\ p\geq 2,
    \label{eq:defdis1}\\
    \Re {\mathscr L}(|u|^{p'-2}u, u) \geq 0 \qquad
    \hbox{\rm if}\ 1<p< 2
      \label{eq:defdis2}
\end{eqnarray}
for all $u\in (\Cspt^{1}(\O))^{m}$. 
Unless otherwise stated we assume that the functions are
complex vector valued.

Saying the $L^{p}$-dissipativity of the operator $A$,
 we mean
the $L^{p}$-dissipa\-tivity of the corresponding form $\elle$,
just to simplify the terminology.

The problem of the dissipativity of linear differential operators and 
 the problem of the contractivity of semigroups 
generated by them attracted
 much attention (see, e.g., 
\cite{mazyasobolev, brezis, davies, amann, strichartz2,
davies1, kovalenko, robinson, davies2, davies3,
liskevich1, liskevichsemen,  langermazya, langer,
daners, karrmann, sobol, liskevich2, ouha, metafune}).
A detailed account of the subject can be found
in the book \cite{ouhabook}, which contains also  
an extensive bibliography.

The present paper is devoted to the  $L^{p}$-dissipativity
($1<p<\infty$)
for partial differential operators. 
It is well known that scalar second order
elliptic operators with
real coefficients may generate contractive semigroups in $L^{p}$
(see \cite{mazyasobolev}).
The case $p=\infty$ was considered in \cite{kresin}, where
necessary and sufficient conditions for the $L^{\infty}$-contractivity
for
scalar second order strongly elliptic systems 
with smooth coefficients were given.
Necessary and sufficient conditions for the
$L^{\infty}$-contractivity 
were later given in \cite{auscher}
under the assumption that the coefficients are measurable and bounded.

The Dirichlet problem for the scalar operator \reff{eq:A}
($m=1$) is considered in \cite{cialmaz} under the assumption that
the entries of $\A$ are complex measures and 
$\Im\A$ is  symmetric.  It is proved that the condition
\begin{equation}
    |p-2| \, |\lan \Im\A \xi,\xi\ran| \leq 2\sqrt{p-1}\lan 
    \Re \A\xi,\xi\ran \quad\quad \forall\ \xi\in\R^{n}
    \label{oldcond}
\end{equation}
 is necessary and sufficient for the 
$L^{p}$-dissipativity.

The main results of the present work are as follows.
In Section \ref{sec:angle}
we use \reff{oldcond} to
obtain the sharp angle of dissipativity of
a scalar complex differential operator $A$. 
To be more precise, we prove in Theorem \ref{th:angdiss} that 
$zA$ ($z\in\C$) is $L^{p}$-dissipative if and only
if $\th_{-}\leq \arg z\leq \th_{+}$, where
$\th_{-}$ and $\th_{+}$ are explicitly given
(see \reff{exprth}). Previously this result was known
for operators with real coefficients
(see \cite{okazawa} and Remark \ref{rem:1} below).
It is worthwhile to remark that we never require 
ellipticity and we may deal with degenerate
matrices. 

In Section \ref{sec:elast},
the two-dimensional elasticity system is considered:
$$
E u = \D u + (1-2\n)^{-1}\nabla\dive u .
$$

After proving a lemma concerning the $L^{p}$-dissipativity
for general systems,
it is shown that $E$ is $L^{p}$-dissipative if and only if
$$ 
\left({1\over 2}-{1\over p}\right)^{2} \leq {2(\n-1)(2\n-1)\over
       (3-4\n)^{2}}
\ .$$

In Section \ref{sec:PDE} we deal with the class
of systems of partial differential equations of the form
$$
    Au=\de_{h}(\A^{h}(x)\de_{h}u)
$$
where $\A^{h}$ are $m\times m$ matrices whose elements
are $L^{1}_{{\rm loc}}$ functions.
We remark that the elasticity system is not of this form.

We find that the operator $A$ is $L^{p}$-dissipative
if and only if
$$\display{
     \Re \lan \A^{h}(x) \l,\l\ran-(1-2/p)^{2}\Re\lan 
     \A^{h}(x)\o,\o\ran (\Re \lan\l,\o\ran)^{2}
     \cr
     -
     (1-2/p)\Re(\lan \A^{h}(x)\o,\l\ran -\lan \A^{h}(x)\l,\o\ran)
     \Re \lan \l,\o\ran   \geq 0
}
$$
 for almost every $x\in\O$ and for any $\l,\o\in\C^{m}$, $|\o|=1$,
 $h=1,\ldots,n$.
 We determine also the angle of dissipativity for such operators.

 In the particular case of  positive real symmetric matrices 
 $\A^{h}$, we
 prove that  $A$ is 
    $L^{p}$-dissipative if and only if
    $$
    \left({1\over 2}-{1\over 
	p}\right)^{2} (\m_{1}^{h}(x)
	+\m_{m}^{h}(x))^{2} \leq \m_{1}^{h}(x)\m_{m}^{h}(x)
    $$
    almost everywhere, $h=1,\ldots,n$,
    where $\m_{1}^{h}(x)$ and $\m_{m}^{h}(x)$ are 
    the smallest and the largest 
    eigenvalues of the matrix $\A^{h}(x)$ respectively.

    The results obtained in Section
\ref{sec:PDE}
are new even for systems of ordinary differential
equations.

\paragrafo{The angle of dissipativity of Second Order 
Scalar Complex Differential Operators} \label{sec:angle}

In this section we consider  the operator
\begin{equation}
    A=\nabla^{t}(\A(x)\nabla)
    \label{cialmazop}
\end{equation}
where $\A =\{a_{ij}(x)\}$ ($i,j=1,\ldots,n$) is a matrix with complex
locally integrable entries defined in a domain $\O\subset \R^{n}$. 
In \cite{cialmaz} it is proved that, if $\Im\A$ is symmetric, there is the 
$L^{p}$-dissipativity of the Dirichlet problem for the differential 
operator $A$ if and only if
\begin{equation}
    |p-2|\, |\lan \Im\A(x)\xi,\xi\ran| \leq 2\sqrt{p-1}\lan 
    \Re\A(x)\xi,\xi\ran
    \label{cialmaz}
\end{equation}
for almost every $x\in\O$ and for any $\xi\in\R^{n}$.

For the sake of completeness we give a proof of the following
elementary  lemma
\begin{lemma}\label{lemma:elem}
    Let $P$ and $Q$ two real measurable functions defined on a set 
    $\O\subset \R^{n}$. 
    Let us suppose that $P(x)\geq 0$ almost everywhere. 
    The inequality
    \begin{equation}
        P(x)\, \cos \th -Q(x)\, \sin\th \geq 0
	\qquad (\th\in [-\pi,\pi])
        \label{diselemen}
    \end{equation}
    holds for almost every $x\in\O$ if and only if
\begin{equation}
        \arccot \, [\essinf_{x\in\Xi} \, (Q(x)/P(x))]-\pi \leq \th \leq
         \arccot \, [\esssup_{x\in\Xi} \, (Q(x)/P(x))]
    \label{tesiel}
\end{equation}
   where $\Xi=\{x\in\O\ |\ P^{2}(x)+Q^{2}(x)>0\}$ and we set
    $$
    Q(x)/P(x)=\cases{
    +\infty  & if $P(x)=0$, $Q(x)>0$\cr
    -\infty  & if $P(x)=0$, $Q(x)<0$.
    }
    $$
    Here  $0<\arccot y <\pi$, $\arccot(+\infty)=0$, 
	$\arccot(-\infty)=\pi$ and 
	$$\essinf_{x\in \Xi} \, (Q(x)/P(x))=+\infty, \quad
\esssup_{x\in \Xi}\, (Q(x)/P(x))=-\infty$$ if $\Xi$ has zero measure.
\end{lemma}

\begin{proof}
    If $\Xi$ has positive measure and $P(x)>0$, inequality \reff{diselemen} means
    $$
    \cos\th - (Q(x)/P(x))\, \sin\th \geq 0
    $$
    and this is true if and only if
    \begin{equation}
	\arccot \, (Q(x)/P(x)) -\pi \leq \th \leq \arccot\, 
	(Q(x)/P(x)).
	\label{condxiel}
    \end{equation}

    If $x\in\Xi$ and $P(x)=0$, \reff{diselemen} means
    $$
    -\pi\leq\th \leq 0, \ \hbox{\rm if}\ Q(x)>0, \quad
    0\leq \th \leq \pi, \ \hbox{\rm if}\ Q(x)<0.
    $$
    This shows that  \reff{diselemen}  is equivalent to \reff{condxiel} 
    provided that $x\in\Xi$. On the other hand, if $x\notin\Xi$,
    $P(x)=Q(x)=0$ almost 
	everywhere and
    \reff{diselemen} is always satisfied. 
    Therefore, if $\Xi$ has positive measure, \reff{diselemen} and
    \reff{tesiel} are equivalent. 

    If $\Xi$ has zero measure, the result is trivial.
\end{proof}

The next Theorem provides a necessary and sufficient condition for the 
$L^{p}$-dissipativity of the Dirichlet problem for the differential 
operator $zA$, where $z\in\C$.

\begin{theorem}\label{th:angdiss}
    Let the matrix $\A$ be symmetric. Let us suppose that the operator $A$ is 
    $L^{p}$-dissipative. Set
    $$\L_{1}=\essinf_{(x,\xi)\in 
	    \Xi} {\lan\Im\A(x)\xi,\xi\ran 
	    \over \lan\Re\A(x)\xi,\xi\ran}, \quad 
	    \L_{2}=\esssup_{(x,\xi)\in 
	    \Xi} {\lan\Im\A(x)\xi,\xi\ran 
	    \over \lan\Re\A(x)\xi,\xi\ran}
	    $$
	    where 
\begin{equation}
            \Xi=\{ (x,\xi)\in \O\times\R^{n}\ |\
	    \lan \Re\A(x)\xi,\xi\ran
            >0\}.
    \label{defXi}
\end{equation}
	    The operator $zA$ is 
    $L^{p}$-dissipative if and only if
   \begin{equation}
        \th_{-}\leq \arg z \leq \th_{+}\, ,
       \label{thcialmaz}
   \end{equation}
       where
       \beginrighe{exprth}
       \th_{-}=\cases{\arccot\left( {2\sqrt{p-1}\over |p-2|} -{p^{2}\over |p-2|}
	       \, {1\over 2\sqrt{p-1}+|p-2|\L_{1}}\right) -\pi
	       & if $p\neq 2$\cr
	       \arccot(\L_{1})-\pi & if $p=2$}
	       \cr
       \th_{+}=\cases{
       \arccot \left(-{2\sqrt{p-1}\over |p-2|} +{p^{2}\over |p-2|}
				      \, 
	{1\over 2\sqrt{p-1}-|p-2|\L_{2}}\right) & if $p\neq 2$\cr
	       \arccot(\L_{2}) & if $p=2$.}
 \endrighe
\end{theorem}

\begin{proof}
    The matrix $\A$ being symmetric, $\Im(e^{i\th}A)$ is symmetric and
    in view of \reff{cialmaz}, the operator $e^{i\th}A$ 
    (with $\th\in[-\pi,\pi]$) is $L^{p}$-dissipative 
    if and only if
\beginrighe{condsincos}
    |p-2|\, |\lan \Re\A(x)\xi,\xi\ran\sin\th+
    \lan \Im\A(x)\xi,\xi\ran \cos\th | \leq \cr
    2\sqrt{p-1}
    (\lan \Re\A(x)\xi,\xi\ran \cos\th -  \lan 
    \Im\A(x)\xi,\xi\ran\sin\th)
\endrighe
    for almost every $x\in\O$ and for any $\xi\in\R^{n}$. Suppose $p\neq 2$.
    Setting
    $$\display{
    a(x,\xi)=|p-2|\, \lan\Re\A(x)\xi,\xi\ran, \quad
    b(x,\xi)=|p-2|\, \lan\Im\A(x)\xi,\xi\ran,\cr
    c(x,\xi)=2\sqrt{p-1}\, \lan\Re\A(x)\xi,\xi\ran, \quad
    d(x,\xi)=2\sqrt{p-1}\, \lan\Im\A(x)\xi,\xi\ran,
    }
    $$
    the 
    inequality in \reff{condsincos} can be written as
    the system
\begin{equation}
    \cases{
      (c(x,\xi)-b(x,\xi))\cos\th-(a(x,\xi)+d(x,\xi)) \sin\th \geq 0, \cr
          (c(x,\xi)+b(x,\xi))\cos\th + (a(x,\xi)-d(x,\xi)) \sin\th\geq 0
          .
    }
    \label{sistemaccio}
\end{equation}

Noting that $c(x,\xi) \pm b(x,\xi)\geq 0$ because of \reff{cialmaz},
the solutions of the inequalities in \reff{sistemaccio} are 
given by the $\th$'s satisfying both of the following conditions
(see Lemma \ref{lemma:elem})
\begin{equation}    
    \cases{\grande \arccot \left(\essinf_{(x,\xi)\in 
    \Xi_{1}}\textstyle{a(x,\xi)+d(x,\xi) \over
    c(x,\xi)-b(x,\xi)}\right)- \pi \leq \th
            \leq \arccot \left(\esssup_{(x,\xi)\in 
    \Xi_{1}}
	    \textstyle{a(x,\xi)+d(x,\xi) \over
    c(x,\xi)-b(x,\xi)}\right)
    \cr\grande
   \arccot \left(\essinf_{(x,\xi)\in 
   \Xi_{2}} \textstyle{d(x,\xi)-a(x,\xi) \over
    c(x,\xi)+b(x,\xi)}\right)- \pi \leq \th
            \leq \arccot \left(\esssup_{(x,\xi)\in 
	    \Xi_{2}} \textstyle{d(x,\xi)-a(x,\xi) \over
    c(x,\xi)+b(x,\xi)}\right) ,
    }
    \label{soluzdis}
\end{equation}
where
$$\display{
\Xi_{1}=\{ (x,\xi)\in \O\times\R^{n}\ |\ (a(x,\xi)+d(x,\xi))^{2}+
(c(x,\xi)-b(x,\xi))^{2}>0\},\cr
\Xi_{2}=\{ (x,\xi)\in \O\times\R^{n}\ |\ (a(x,\xi)-d(x,\xi))^{2}+
(b(x,\xi)+c(x,\xi))^{2}>0\}.
}
$$

 We have
$$\display{
a(x,\xi)\, d(x,\xi)= b(x,\xi)\, c(x,\xi),\cr
a^{2}(x,\xi)+b^{2}(x,\xi)+c^{2}(x,\xi)+d^{2}(x,\xi)=
p^{2}(\lan \Re\A(x)\xi,\xi\ran^{2}+
    \lan \Im\A(x)\xi,\xi\ran^{2})
}
$$
and then, keeping in mind \reff{cialmaz},
we may write $\Xi_{1}=\Xi_{2}=\Xi$, where
$\Xi$ is given by \reff{defXi}.

Moreover
$$
{a(x,\xi)+d(x,\xi) \over
    c(x,\xi)-b(x,\xi)} \geq {d(x,\xi)-a(x,\xi) \over
    c(x,\xi)+b(x,\xi)}
$$
and  then $\th$ satisfies all of the inequalities in \reff{soluzdis} if 
and only if
\begin{equation}
    \arccot \left(\essinf_{(x,\xi)\in 
       \Xi} \textstyle{d(x,\xi)-a(x,\xi) \over
        c(x,\xi)+b(x,\xi)}\right)- \pi \leq
        \th \leq
        \arccot \left(\esssup_{(x,\xi)\in 
        \Xi}\textstyle{a(x,\xi)+d(x,\xi) \over
        c(x,\xi)-b(x,\xi)}\right)
    \label{primacondsup}
\end{equation}

A direct computation shows that
$$\display{
{d(x,\xi)-a(x,\xi) \over
	c(x,\xi)+b(x,\xi)}= {2\sqrt{p-1}\over |p-2|} -{p^{2}\over |p-2|}
	\, {1\over 2\sqrt{p-1}+|p-2|\L(x,\xi)},\cr
	{a(x,\xi)+d(x,\xi) \over
		c(x,\xi)-b(x,\xi)} = 
		-{2\sqrt{p-1}\over |p-2|} +{p^{2}\over |p-2|}
			\, {1\over 2\sqrt{p-1}-|p-2|\L(x,\xi)}
	}
$$
where
$$
\L(x,\xi)= {\lan\Im\A(x)\xi,\xi\ran \over \lan\Re\A(x)\xi,\xi\ran}\, .
$$

Hence condition \reff{primacondsup} is satisfied if and only if
\reff{thcialmaz} holds.

If $p=2$, \reff{condsincos} is simply
$$ 
    \lan \Re\A(x)\xi,\xi\ran \cos\th -  \lan 
    \Im\A(x)\xi,\xi\ran\sin\th \geq 0
$$
and the result follows directly from Lemma \ref{lemma:elem}.
\end{proof}

\begin{remark}\label{rem:1} If $\A$ is a real matrix, then $\L_{1}=\L_{2}=0$ and the
angle of dissipativity does not depend on the operator. In fact 
we have
$$
{2\sqrt{p-1}\over |p-2|} -{p^{2}\over 2\sqrt{p-1}|p-2|}=
- {|p-2|\over 2\sqrt{p-1}}
$$
and Theorem \ref{th:angdiss} shows that $zA$ 
is dissipative if and only 
if
$$
\arccot \left(- {|p-2|\over 2\sqrt{p-1}}\right) -\pi \leq
\arg z \leq \arccot\left({|p-2|\over 2\sqrt{p-1}}\right),
$$
i.e.
$$
|\arg z| \leq \arctan \left({2\sqrt{p-1}\over |p-2|}\right).
$$

This is a well known result (see, e.g., \cite{fattorini}, 
\cite{fattorini2},
\cite{okazawa}).
\end{remark}

\paragrafo{Two-dimensional Elasticity} \label{sec:elast}

Let us consider the classical operator of 
two-dimensional elasticity
\begin{equation}
    Eu=\D u + (1-2\n)^{-1}\nabla \dive u
    \label{opelast}
\end{equation}
where  $\n$ is the Poisson ratio. 
It is well known that $E$ 
 is strongly elliptic if and only if either $\n>1$ or $\n<1/2$. 

In this Section we give a necessary and sufficient condition for 
the $L^{p}$-dissipativity of operator \reff{opelast}.

We start giving a necessary condition for the 
$L^{p}$-dissipativity of the operator
\begin{equation}
    A=\de_{h} (\A^{hk}(x)\de_{k})
    \label{eq:A2}
\end{equation}
where  $\A^{hk}(x)=\{a^{hk}_{ij}(x)\}$ are 
$m\times m$ matrices 
whose elements are complex locally integrable  functions defined in 
an arbitrary domain $\O$ of $\R^{2}$ 
$(1\leq i,j\leq 
m,\ 1\leq h,k \leq 2)$. 

The following lemma holds in any number of variables.

\begin{lemma}\label{lemma:5}
    Let $\O$ be a domain of $\R^{n}$.
    The operator \reff{eq:A2} is $L^{p}$-dissipative if and only if
    \beginrighe{cond1}
    \int_{\O}\Big(\Re \lan \A^{hk} \de_{k}v,\de_{h}v\ran
    \cr
    -(1-2/p)^{2}|v|^{-4}\Re\lan \A^{hk} v,v\ran \Re \lan 
    v,\de_{k}v\ran \Re \lan v,\de_{h}v\ran \cr
    -(1-2/p)|v|^{-2}\Re(\lan \A^{hk} v,\de_{h}v\ran \Re \lan v,\de_{k}v\ran
   \cr
    -\lan \A^{hk} \de_{k}v,v\ran \Re \lan v,\de_{h}v\ran)
     \Big) dx 
    \geq 0
   \endrighe
    for any $v\in (\Cspt^{1}(\O))^{m}$. Here and in the sequel the 
      integrand is extended by zero on the set where $v$ vanishes.
\end{lemma}

\begin{proof}
    \textit{Sufficiency.} First suppose $p\geq 2$. 
       Let $u\in (\Cspt^{1}(\O))^{m}$ and set 
       $v=|u|^{(p-2)/2}u$. We have $v\in (\Cspt^{1}(\O))^{m}$ and
       $u=|v|^{(2-p)/p}v$. From the identities
       $$\display{
       \lan \A^{hk} \de_{k}u, \de_{h}(|u|^{p-2}u)\ran =\cr
       \lan \A^{hk} \de_{k}v,\de_{h}v\ran 
	   -(1-2/p)^{2}|v|^{-2}\Re\lan \A^{hk} v,v\ran\,
	  \de_{k}|v|\de_{h}|v| \cr
	   -(1-2/p)|v|^{-1}\Re(\lan \A^{hk} v,\de_{h}v\ran\, 
	   \de_{k}|v|
	   -\lan \A^{hk} \de_{k}v,v\ran \de_{h}|v|),
	   \cr
	   }
       $$
       $$
       \de_{k}|v|= |v|^{-1}\Re\lan v,\de_{k}v\ran,
       $$
       we see that the left-hand side in \reff{cond1} is equal
       to $\elle(u,|u|^{p-2}u)$. Then \reff{eq:defdis1} is satisfied for 
       any $u\in (\Cspt^{1}(\O))^{m}$.
       
       If $1<p<2$, we may write \reff{eq:defdis2} as
       $$\display{
       \Re \int_{\O}\lan (\A^{hk})^{*}\de_{h}u, 
       \de_{k}(|u|^{p'-2}u)\ran \, dx \geq 0
       }
       $$
       for any $u\in(\Cspt^{1}(\O))^{m}$. The first part of the proof shows 
       that this implies
       \beginrighe{cond1*}
       \int_{\O}\Big(\Re \lan (\A^{hk})^{*} \de_{h}v,\de_{k}v\ran 
   \cr
       -(1-2/p')^{2}|v|^{-4}\Re\lan (\A^{hk})^{*} v,v\ran \Re \lan 
       v,\de_{h}v\ran \Re \lan v,\de_{k}v\ran \cr
       -(1-2/p')|v|^{-2}\Re(\lan (\A^{hk})^{*} v,\de_{k}v\ran 
       \Re \lan v,\de_{h}v\ran
\cr
       -\lan (\A^{hk})^{*} \de_{h}v,v\ran \Re \lan v,\de_{k}v\ran)
	\Big) dx 
       \geq 0
      \endrighe
       for any $v\in(\Cspt^{1}(\O))^{m}$. 
     Since $1-2/p'=-(1-2/p)$, this inequality is exactly 
     \reff{cond1}.
     
     \textit{Necessity.} Let $p\geq 2$ and set
      $$
      g_{\e}=(|v|^{2}+\e^{2})^{1/2}, \quad u_{\e}=g_{\e}^{2/p-1}v,
      $$
      where  $v\in(\Cspt^{1}(\O))^{m}$.  We have
      $$\display{
      \lan \A^{hk}\de_{k} u_{\e}, \de_{h}(|u_{\e}|^{p-2}u_{\e})\ran 
      \cr
      =|u_{\e}|^{p-2}\lan \A^{hk}\de_{k} u_{\e}, \de_{h}u_{\e}\ran +
      (p-2)|u_{\e}|^{p-3}\lan \A^{hk}\de_{k} u_{\e},u_{\e}\ran 
      \de_{h}|u_{\e}|.
      }
      $$
      
      One checks directly that
      $$\display{
      |u_{\e}|^{p-2}\lan \A^{hk}\de_{k} u_{\e}, \de_{h}u_{\e}\ran
      \cr
      =
      (1-2/p)^{2}g_{\e}^{-(p+2)}|v|^{p-2}\lan \A^{hk} v,v\ran \, 
      \Re\lan 
	      v,\de_{k}v\ran\, \Re\lan 
	      v,\de_{h}v\ran\cr
      -(1-2/p)g_{\e}^{-p}|v|^{p-2}(\lan \A^{hk} v,\de_{h}v\ran 
      \Re\lan 
	v,\de_{k}v\ran+\lan\A^{hk}\de_{k}v,v\ran
	\Re\lan 
	v,\de_{h}v\ran)
      \cr
      + g_{\e}^{2-p}|v|^{p-2}\lan \A^{hk}\de_{k} v,\de_{h}v\ran ,
      \cr
      |u_{\e}|^{p-3}\lan \A^{hk}\de_{k} u_{\e},u_{\e}\ran 
	   \, \de_{h}|u_{\e}|
      \cr
     = (1-2/p)[(1-2/p)g_{\e}^{-(p+2)}|v|^{p-2}
     \cr - 
     g_{\e}^{-p}|v|^{p-4}]\, \lan \A^{hk} 
      v,v\ran \, \Re\lan 
	v,\de_{k}v\ran\, \Re\lan 
	v,\de_{h}v\ran\cr
      + [g_{\e}^{2-p}|v|^{p-4}-(1-2/p)g_{\e}^{-p}|v|^{p-2}]\, 
      \lan \A^{hk}\de_{k}v,v\ran\, \Re\lan 
      v,\de_{h}v\ran
      }
      $$
      on the set $E=\{x\in \O\ |\ |v(x)|>0\}$.
      The inequality $g_{\e}^{a}\leq |v|^{a}$ for $a\leq 0$, shows that
	the right-hand sides are majorized by $L^{1}$ functions.
      Since $g_{\e}\to |v|$ pointwise as $\e\to 0^{+}$, we find
      $$\display{
      \lim_{\e\to 0^{+}}  \lan \A^{hk}\de_{k} u_{\e}, 
      \de_{h}(|u_{\e}|^{p-2}u_{\e})\ran
      \cr
      =\lan \A^{hk} \de_{k}v,\de_{h}v\ran 
    -(1-2/p)^{2}|v|^{-4}\lan \A^{hk} v,v\ran \Re \lan 
    v,\de_{k}v\ran \Re \lan v,\de_{h}v\ran \cr
    -(1-2/p)|v|^{-2}(\lan \A^{hk} v,\de_{h}v\ran \Re \lan v,\de_{k}v\ran
    -\lan \A^{hk} \de_{k}v,v\ran \Re \lan v,\de_{h}v\ran)
     }
     $$
      and dominated convergence gives
      \beginrighe{wefind}
      \lim_{\e\to 0^{+}} \elle(u_{\e},|u_{\e}|^{p-2}u_{\e}) = 
      \lim_{\e\to 0^{+}} \int_{E}
      \lan \A^{hk}\de_{k} u_{\e}, 
	    \de_{h}(|u_{\e}|^{p-2}u_{\e})\ran
     \, dx =
      \cr\grande
      \Re\int_{E}[\lan \A^{hk} \de_{k}v,\de_{h}v\ran
     \cr
    -(1-2/p)^{2}|v|^{-4}\lan \A^{hk} v,v\ran \Re \lan 
    v,\de_{k}v\ran \Re \lan v,\de_{h}v\ran \cr
    -(1-2/p)|v|^{-2}(\lan \A^{hk} v,\de_{h}v\ran \Re \lan v,\de_{k}v\ran
    \cr
    -\lan \A^{hk} \de_{k}v,v\ran \Re \lan v,\de_{h}v\ran)]\, dx.
     \endrighe
     
      The function $u_{\e}$ being in $(\Cspt^{1}(\O))^{m}$, 
      \reff{eq:defdis1} implies \reff{cond1}.
      
      If $1<p<2$, from \reff{wefind} it follows that
      $$\display{
      \lim_{\e\to 0^{+}} \elle(|u_{\e}|^{p'-2}u_{\e},u_{\e}) = 
      \lim_{\e\to 0^{+}}\, \Re \int_{E} \lan (\A^{hk})^{*}\de_{h}u, 
       \de_{k}(|u|^{p'-2}u)\ran \, dx =
       \cr
       \int_{E}\Big(\Re \lan (\A^{hk})^{*} \de_{h}v,\de_{k}v\ran 
    \cr
	-(1-2/p')^{2}|v|^{-4}\Re\lan (\A^{hk})^{*} v,v\ran \Re \lan 
	v,\de_{h}v\ran \Re \lan v,\de_{k}v\ran \cr
	-(1-2/p')|v|^{-2}\Re(\lan (\A^{hk})^{*} v,\de_{k}v\ran 
	\Re \lan v,\de_{h}v\ran
 \cr
	-\lan (\A^{hk})^{*} \de_{h}v,v\ran \Re \lan v,\de_{k}v\ran)
	 \Big) dx.
       }
       $$
      
       This shows that \reff{eq:defdis2} implies \reff{cond1*} and
       the proof is complete. 
\end{proof}

\begin{theorem}\label{th:5}
    Let $\O$ be a domain of $\R^{2}$.
    If the operator \reff{eq:A2} is $L^{p}$-dissipative, we have 
    \beginrighe{algcond}
   \Re  \lan (\A^{hk}(x)\xi_{h}\xi_{k})\l,\l\ran  -(1-2/p)^{2}\Re\lan (\A^{hk}(x)\xi_{h}\xi_{k})\o,\o\ran (\Re \lan 
    \l,\o\ran)^{2}\cr
    -(1-2/p)\Re(\lan (\A^{hk}(x)\xi_{h}\xi_{k})\o,\l\ran
	-\lan (\A^{hk}(x)\xi_{h}\xi_{k})\l,\o\ran) 
	\Re \lan \l,\o\ran 
	\cr\geq 0
	\endrighe
	for almost every $x\in\O$ and 
	for any $\xi\in\R^{2}$, $\l,\, \o\in \C^{m}$, $|\o|=1$.
\end{theorem}

\begin{proof}
   Let us assume that $\A$ is a constant matrix and that $\O=\R^{2}$.
   Let us fix $\o\in\C^{m}$ with $|\o|=1$ and take
   $
   v(x)=w(x)\, \h(\log|x|/\log R)
   $,
   where 
   \begin{equation}
       w(x)=\m\o+\p(x),
       \label{defw}
   \end{equation}
   $\m,\,R\in\R^{+}$, $R>1$,
   $\p\in(\Cspt^{\infty}(\R^{2}))^{m}$, 
   $\h\in C^{\infty}(\R)$, $\h(t)=1$ if $t\leq 1/2$ and
   $\h(t)=0$ if $t\geq 1$.
   
   We have
   $$\display{
   \lan \A^{hk}\de_{k}v,\de_{h}v\ran = \lan 
   \A^{hk}\de_{k}w,\de_{h}w\ran \h^{2}(\log |x|/\log R)\cr
   + (\log R)^{-1}(\lan \A^{hk}\de_{k}w,w\ran x_{h} +
   \lan \A^{hk}w,\de_{h}w\ran x_{k})\times\cr
   |x|^{-2}\h(\log |x|/\log R)\,
   \h'(\log |x|/\log R)\cr
   + (\log R)^{-2}\lan \A^{hk}w,w\ran x_{h}x_{k}|x|^{-4}
   \left(\h'(\log |x|/\log R)\right)^{2}
   }
   $$
   and then, choosing $\d$ such that $\spt \p\subset B_{\d}(0)$,
   $$\display{
   \int_{\R^{2}} \lan \A^{hk}\de_{k}v,\de_{h}v\ran dx =
   \int_{B_{\d}(0)}\lan \A^{hk}\de_{k}w,\de_{h}w\ran dx\cr
   +
   {1\over \log^{2}R}\int_{B_{R}(0)\setminus B_{\sqrt{R}}(0)} 
   \lan \A^{hk}w,w \ran {x_{h}x_{k}\over |x|^{4}}
      \left(\h'(\log |x|/\log R)\right)^{2}dx
   }
   $$
   provided that $R>\d^{2}$.
   Since
   $$
   \lim_{R\to +\infty}{1\over \log^{2}R}\int_{B_{R}(0)\setminus B_{\sqrt{R}}(0)} 
   {dx\over |x|^{2}}
       = 0,
   $$
   we have
   $$
   \lim_{R\to +\infty} \int_{\R^{2}} \lan \A^{hk}\de_{k}v,\de_{h}v\ran dx =
   \int_{B_{\d}(0)}\lan \A^{hk}\de_{k}w,\de_{h}w \ran dx.
   $$
   
   On the set where $v\neq 0$ we have
   $$\display{
   |v|^{-4}\lan\A^{hk}v,v\ran \Re \lan v,\de_{k}v\ran \Re \lan 
   v,\de_{h}v\ran =\cr
   |w|^{-4}\lan\A^{hk}w,w\ran \Re \lan w,\de_{k}w\ran \Re \lan 
   w,\de_{h}w\ran \h^{2}(\log|x|/ \log R) \cr
   + (\log R)^{-1}|w|^{-2}\lan\A^{hk}w,w\ran  (
   \Re \lan w,\de_{h}w\ran x_{k}+ 
   \Re \lan w,\de_{k}w\ran x_{h})|x|^{-2}\times\cr
   \h(\log|x|/ \log R) \,
   \h'(\log|x|/ \log R) \cr
   + (\log R)^{-2}\lan\A^{hk}w,w\ran x_{h}x_{k}|x|^{-4}
   (\h'(\log|x|/\log R))^{2}
   }
   $$
   and then
   $$\display{
  \lim_{R\to+\infty}\int_{\R^{2}}|v|^{-4}\lan \A^{hk}v,v\ran \Re \lan 
    v,\de_{k}v\ran \Re \lan v,\de_{h}v\ran dx= \cr
   \int_{B_{\d}(0)}
   |w|^{-4}\lan\A^{hk}w,w\ran \Re \lan w,\de_{k}w\ran \Re \lan 
      w,\de_{h}w\ran dx .
   }
   $$
   
   In the same way we obtain
   $$\display{
   \lim_{R\to+\infty}\int_{\R^{2}}|v|^{-2}\Re(\lan \A^{hk}v,\de_{h}v\ran \Re \lan v,\de_{k}v\ran
    -\lan \A^{hk}\de_{k}v,v\ran \Re \lan v,\de_{h}v\ran)dx=\cr
    \int_{B_{\d}(0)}|w|^{-2}\Re(\lan \A^{hk}w,\de_{h}w\ran \Re \lan w,\de_{k}w\ran
	-\lan \A^{hk}\de_{k}w,w\ran \Re \lan w,\de_{h}w\ran)dx.
   }
   $$
   
   In view of Lemma \ref{lemma:5}, \reff{cond1} holds. 
   Putting $v$ in this formula and letting $R\to+\infty$, we find
   \beginrighe{cond1w}
       \int_{B_{\d}(0)}\Big(\Re \lan \A^{hk}\de_{k}w,\de_{h}w\ran \cr
       -(1-2/p)^{2}|w|^{-4}\Re\lan \A^{hk}w,w\ran \Re \lan 
       w,\de_{k}w\ran \Re \lan w,\de_{h}w\ran \cr
       -(1-2/p)|w|^{-2}\Re(\lan \A^{hk}w,\de_{h}w\ran \Re \lan w,\de_{k}w\ran
       \cr
       -\lan \A^{hk}\de_{k}w,w\ran \Re \lan w,\de_{h}w\ran)
	\Big) dx 
       \geq 0.
      \endrighe
   
   On the other hand, keeping in mind \reff{defw}, 
   $$\display{\Re \lan \A^{hk}\de_{k}w,\de_{h}w\ran =
   \Re \lan \A^{hk}\de_{k}\p,\de_{h}\p \ran ,\cr 
      |w|^{-4}\Re\lan \A^{hk}w,w\ran \Re \lan 
       w,\de_{k}w\ran \Re \lan w,\de_{h}w\ran =\cr
       |\m\o+\p|^{-4}\Re\lan \A^{hk}(\m\o+\p),\m\o+\p\ran \Re \lan 
	      \m\o+\p,\de_{k}\p\ran \Re \lan \m\o+\p,\de_{h}\p\ran, \cr
     |w|^{-2}\Re(\lan \A^{hk}w,\de_{h}w\ran \Re \lan w,\de_{k}w\ran 
     -\lan \A^{hk}\de_{k}w,w\ran \Re \lan w,\de_{h}w\ran)=\cr
     |\m\o+\p|^{-2}\Re(\lan \A^{hk}(\m\o+\p),\de_{h}\p\ran \Re \lan 
     \m\o+\p,\de_{k}\p\ran \cr
     -\lan \A^{hk}\de_{k}(\m\o+\p),\m\o+\p\ran \Re \lan 
     \m\o+\p,\de_{h}\p\ran).
       \cr  
  }
  $$
  
  Letting $\m\to+\infty$ in \reff{cond1w}, we obtain
  \beginrighe{cond1psi}
  \int_{\R^{2}}\Big(\Re \lan \A^{hk}\de_{k}\p,\de_{h}\p\ran
	  \cr
	  -(1-2/p)^{2}\Re\lan \A^{hk}\o,\o\ran \Re \lan 
	  \o,\de_{k}\p\ran \Re \lan \o,\de_{h}\p\ran \cr
	  -(1-2/p)\Re(\lan \A^{hk}\o,\de_{h}\p\ran \Re \lan \o,\de_{k}\p\ran
  \cr	-\lan \A^{hk}\de_{k}\p,\o\ran \Re \lan \o,\de_{h}\p\ran)
	  \Big)\, dx  \geq 0. 
\endrighe

Putting in \reff{cond1psi}
$$
\p(x)=\l\,\f(x)\, e^{i\m\lan\xi,x\ran}
$$
where $\l\in\C^{m}$, $\f\in\Cspt^{\infty}(\R^{2})$ and $\m$ is a real parameter,
by standard arguments (see, e.g., \cite[p.107--108]{fichera}),  
we find \reff{algcond}.

If the matrix $\A$ is not constant, 
   take $\p \in 
    (\Cspt^{1}(\R^{2}))^{m}$ and define
    $$
    v(x)=\p((x-x_{0})/\e)
    $$
    where $x_{0}$ is a fixed point in $\O$ and $0<\e<\hbox{dist}\, 
    (x_{0},\de\O)$. 
    
    Putting this particular $v$ in \reff{cond1} and making a change 
    of variables, we obtain
    $$\display{
	\int_{\R^{2}}\Big(\Re \lan \A^{hk}(x_0+\e y)\de_{k}\p,\de_{h}\p\ran 
	\cr
	-(1-2/p)^{2}|\p|^{-4}\Re\lan \A^{hk}(x_0+\e y)\p,\p\ran \Re \lan 
	\p,\de_{k}\p\ran \Re \lan \p,\de_{h}\p\ran \cr
	-(1-2/p)|\p|^{-2}\Re(\lan \A^{hk}(x_0+\e y)\p,\de_{h}\p\ran \Re \lan \p,\de_{k}\p\ran
\cr	-\lan \A^{hk}(x_0+\e y)\de_{k}\p,\p\ran \Re \lan \p,\de_{h}\p\ran)
	\Big)\, dy  \geq 0.
}      
	$$
    
	Letting $\e \to 0^{+}$ we find
	$$\display{
	       \int_{\R^{2}}\Big(\Re \lan 
	       \A^{hk}(x_0)\de_{k}\p,\de_{h}\p\ran \cr
	       -(1-2/p)^{2}|\p|^{-4}\Re\lan \A^{hk}(x_0)\p,\p\ran \Re \lan 
	       \p,\de_{k}\p\ran \Re \lan \p,\de_{h}\p\ran \cr
	       -(1-2/p)|\p|^{-2}\Re(\lan \A^{hk}(x_0)\p,\de_{h}\p\ran \Re \lan \p,\de_{k}\p\ran
       \cr	-\lan \A^{hk}(x_0)\de_{k}\p,\p\ran \Re \lan \p,\de_{h}\p\ran)
	       \Big)\, dy  \geq 0
 }
 $$
		for almost every $x_{0}\in\O$. The arbitrariness of $\p\in 
		(\Cspt^{1}(\R^{2}))^{m}$ and what we have proved for constant 
		matrices give the result.
\end{proof}

Since in problem of Elasticity we are  interested
in real solutions, we shall discuss the $L^{p}$-dissipativity of
the operator \reff{opelast} in a real frame.
In the present Section, all the functions we are going to consider,
in particular the ones appearing  
in the conditions \reff{eq:defdis1} and \reff{eq:defdis2},
are supposed to be real vector valued.

\begin{theorem}\label{th:cnesdiss}
   The operator \reff{opelast} is $L^{p}$-dissipative if
   and only if 
  \begin{equation}
      \left({1\over 2}-{1\over p}\right)^{2} \leq 
      {2(\n-1)(2\n-1)\over
	    (3-4\n)^{2}}. 
       \label{condvecchia}
   \end{equation}
\end{theorem}

\begin{proof}    
    \textit{Necessity.}  We have
    $$\display{
    \lan (\A^{hk}\xi_{h}\xi_{k})\l,\l\ran =
    |\xi|^{2}|\l|^{2}+(1-2\n)^{-1}\lan\xi,\l\ran^{2}, \cr
    \lan (\A^{hk}\xi_{h}\xi_{k})\o,\o\ran =
	|\xi|^{2}+(1-2\n)^{-1}\lan\xi,\o\ran^{2},\cr
	\lan (\A^{hk}\xi_{h}\xi_{k})\l,\o\ran =
	    |\xi|^{2}|\lan \l,\o\ran +(1-2\n)^{-1}\lan\xi,\l\ran
	    \lan\xi,\o\ran
    }
    $$
    for any $\xi,\, \l,\, \o\in\R^{2}$, $|\o|=1$. Hence, in view of 
    Theorem \ref{th:5}, the $L^{p}$-dissipativity of $E$ implies
    \beginrighe{condalgelas}
-(1-2/p)^{2}[|\xi|^{2}+(1-2\n)^{-1}(\xi_{j}\o_{j})^{2}](\l_{j}\o_{j})^{2}
\cr
+
|\xi|^{2}|\l|^{2} + (1-2\n)^{-1}(\xi_{j}\l_{j})^{2}\geq 0
    \endrighe
for any $\xi,\, \l,\, \o\in\R^{2}$, $|\o|=1$.   
    
Without loss of generality 
 we may suppose  $\xi=(1,0)$. Setting
$C_{p}=(1-2/p)^{2}$ and $\g= (1-2\n)^{-1}$, condition \reff{condalgelas} 
can be written as
\begin{equation}
     -C_{p}(1+\g \o_{1}^{2})(\l_{j}\o_{j})^{2}+
     |\l|^{2}+\g \l_{1}^{2}\geq 0
    \label{quadrform}
\end{equation}
 for any $\l,\, \o\in\R^{2}$, $|\o|=1$. 

Condition \reff{quadrform} holds if and only if
$$\display{
-C_{p}(1+\g\o_{1}^{2})\o_{1}^{2}+1+\g\geq 0,
\cr
[C_{p}(1+\g\o_{1}^{2})\o_{1}\o_{2}]^{2}\leq \cr
[-C_{p}(1+\g\o_{1}^{2})\o_{1}^{2}+1+\g]\, 
[-C_{p}(1+\g\o_{1}^{2})\o_{2}^{2}+1]
}
$$
for any $\o\in\R^{2}$, $|\o|=1$. 

In particular, the second condition 
has to be satisfied.
This can be written in the form
\begin{equation}
    1+\g -C_{p}(1+\g\o_{1}^{2})(1+\g\o_{2}^{2})\geq 0
    \label{finalform}
\end{equation}
for any $\o\in\R^{2}$, $|\o|=1$. The minimum of the left hand side of 
\reff{finalform} on 
the unit sphere is given by
$$
1+\g-C_{p}(1+\g/2)^{2}.
$$
Hence \reff{finalform} is satisfied if and only if
$1+\g-C_{p}(1+\g/2)^{2}\geq 0$. The last inequality means
$$
{2(1-\n)\over 1-2\n} - \left({p-2\over p}\right)^{2}\left({3-4\n\over 2(1-2\n)}\right)^{2}
\geq 0
$$
i.e. \reff{condvecchia}.
From the identity $4/(p\,p')=1-(1-2/p)^{2}$, it follows
that \reff{condvecchia} can be written also as
\begin{equation}
     {4\over p\, p'} \geq {1\over (3-4\n)^{2}}. 
      \label{condelast}
  \end{equation}

\textit{Sufficiency.} In view of Lemma
\ref{lemma:5}, $E$ is $L^{p}$-dissipative if and only if
\begin{equation}
    \int_{\O}[-C_{p}|\nabla|v||^{2}+\sum_{j}|\nabla v_{j}|^{2}
    -\g\, C_{p}\, |v|^{-2}|v_{h}\de_{h}|v||^{2}+\g\, |\dive v|^{2}]\, 
    dx \geq 0
    \label{eq:condizv}
\end{equation}
for any $v\in (\Cspt^{1}(\O))^{2}$. 
Choose $v\in (\Cspt^{1}(\O))^{2}$ and define
$$\display{
X_{1}=|v|^{-1}(v_{1}\de_{1}|v|+v_{2}\de_{2}|v|),\quad 
X_{2}=|v|^{-1}(v_{2}\de_{1}|v|-v_{1}\de_{2}|v|) \cr
Y_{1}=|v|[\de_{1}(|v|^{-1}v_{1})+\de_{2}(|v|^{-1}v_{2})],
\quad
Y_{2}=|v|[\de_{1}(|v|^{-1}v_{2})-\de_{2}(|v|^{-1}v_{1})]
}
$$
on the set $E=\{x\in\O\ |\ v\neq 0\}$.
From the identities
$$
|\nabla|v||^{2}=X_{1}^{2}+X_{2}^{2}, \quad
Y_{1}=(\de_{1}v_{1}+\de_{2}v_{2}) - X_{1}, \quad
Y_{2}=(\de_{1}v_{2}-\de_{2}v_{1}) - X_{2}
$$
it follows
$$\display{
Y_{1}^{2}+Y_{2}^{2}=|\nabla|v||^{2}+
(\de_{1}v_{1}+\de_{2}v_{2})^{2}+
(\de_{1}v_{2}-\de_{2}v_{1})^{2}\cr
-2(\de_{1}v_{1}+\de_{2}v_{2})X_{1}-
2(\de_{1}v_{2}-\de_{2}v_{1})X_{2}.
}
$$

Keeping in mind that $\de_{h}|v|=|v|^{-1}v_{j}\de_{h}v_{j}$, one
can check that
$$\display{
(\de_{1}v_{1}+\de_{2}v_{2})(v_{1}\de_{1}|v|+v_{2}\de_{2}|v|)+
(\de_{1}v_{2}-\de_{2}v_{1})(v_{2}\de_{1}|v|-v_{1}\de_{2}|v|)=\cr
|v|\, |\nabla |v||^{2}+ |v|(\de_{1}v_{1}
\de_{2}v_{2}-\de_{2}v_{1}\de_{1}v_{2}),
}
$$
which implies
\begin{equation}
    \sum_{j}|\nabla v_{j}|^{2} = X_{1}^{2}+X_{2}^{2}+
    Y_{1}^{2}+Y_{2}^{2}.
    \label{relXY}
\end{equation}

Thus \reff{eq:condizv} can be written as
\begin{equation}
    \int_{E}\left[{4\over p\, p'}(X_{1}^{2}+X_{2}^{2}) + Y_{1}^{2}+Y_{2}^{2}
    -\g\, C_{p}X_{1}^{2}+\g\, (X_{1}+Y_{1})^{2}\right]dx \geq 0.
    \label{premagic}
\end{equation}

Let us prove that
\begin{equation}
    \int_{E}X_{1}Y_{1}dx = -\int_{E}X_{2}Y_{2}dx.
    \label{magic}
\end{equation}

Since $X_{1}+Y_{1}=\dive v$ and 
$X_{2}+Y_{2}=\de_{1}v_{2}-\de_{2}v_{1}$, keeping in mind 
\reff{relXY}, we may write
$$\display{
2\int_{E}(X_{1}Y_{1}+X_{2}Y_{2})dx =
\int_{E}[(X_{1}+Y_{1})^{2}+(X_{2}+Y_{2})^{2}-
(X_{1}^{2}+X_{2}^{2}+
    Y_{1}^{2}+Y_{2}^{2})]\, dx=\cr
\int_{E}[(\dive v)^{2}+
(\de_{1}v_{2}-\de_{2}v_{1})^{2}-
\sum_{j}|\nabla v_{j}|^{2}]\, dx
}
$$
i.e.
$$
\int_{E}(X_{1}Y_{1}+X_{2}Y_{2})dx =
\int_{E}(
\de_{1}v_{1}\de_{2}v_{2}-
\de_{1}v_{2}\de_{2}v_{1})\,dx.
$$

The set $\{x\in \O\setminus E\ |\ \nabla v(x)\neq 0\}$ has zero 
measure and then 
$$
\int_{E}(X_{1}Y_{1}+X_{2}Y_{2})dx =
\int_{\O}(
\de_{1}v_{1}\de_{2}v_{2}-
\de_{1}v_{2}\de_{2}v_{1})\,dx.
$$

There exists a sequence $\{v^{(n)}\}\subset \Cspt^{\infty}(\O)$ such 
that $v^{(n)}\to v$, $\nabla v^{(n)}\to \nabla v$ uniformly in $\O$ 
and hence
$$\display{
\int_{\O}\de_{1}v_{1}\de_{2}v_{2}dx = \lim_{n\to\infty}
\int_{\O}\de_{1}v^{(n)}_{1}\de_{2}v^{(n)}_{2}dx=\cr
\lim_{n\to\infty}\int_{\O}\de_{1}v^{(n)}_{2}\de_{2}v^{(n)}_{1}dx=
\int_{\O}\de_{1}v_{2}\de_{2}v_{1}dx
}
$$
 and \reff{magic} is proved.
In view of this, \reff{premagic} can be written as
$$\display{
\int_{E}\left({4\over p\, p'}(1+\g)X_{1}^{2}+2\th \g\, X_{1}Y_{1}
+(1+\g)Y_{1}^{2}\right) dx\cr
+
\int_{E}\left({4\over p\, p'}X_{2}^{2} -2(1-\th) \g\, X_{2}Y_{2}
+Y_{2}^{2}\right) dx \geq 0
}
$$
for any fixed $\th\in\R$.

If we choose
$$
\th={2(1-\n)\over 3-4\n}
$$
we find
$$
(1-\th)\g={1\over 3-4\n}\, , \qquad
\th^{2}\g^{2}={(1+\g)^{2}\over (3-4\n)^{2}}\, 
.
$$

Inequality \reff{condelast} leads to
$$
\th^{2}\g^{2}\leq {4\over p\,p'}\, (1+\g)^{2}, \qquad
(1-\th)^{2}\g^{2}\leq {4\over p\,p'}\, .
$$

Observing that \reff{condvecchia} implies $1+\g=
2(1-\n)(1-2\n)^{-1}\geq 0$, we get
$$\display{
{4\over p\, p'}(1+\g)x_{1}^{2}+2\th \g\, x_{1}y_{1}
+(1+\g)y_{1}^{2}\geq 0,\cr
{4\over p\, p'}x_{2}^{2} -2(1-\th) \g\, x_{2}y_{2}
+y_{2}^{2} \geq 0
}
$$
for any $x_{1},x_{2},y_{1},y_{2}\in\R$.  This shows that
\reff{premagic} holds. Then 
\reff{eq:condizv} is true 
 for any $v\in (\Cspt^{1}(\O))^{2}$ and the proof is complete.
\end{proof}

We shall now give two Corollaries of this result. They concerns
the comparison between $E$ and $\D$ from the point of view
of the $L^{p}$-dissipativity.

\begin{corollary}\label{th:7}
  There exists $k>0$ such that $E-k\D$ is $L^{p}$-dissipative if 
    and only if 
    \begin{equation}
	\left({1\over 2}-{1\over p}\right)^{2} < 
	     {2(\n-1)(2\n-1)\over
		   (3-4\n)^{2}}\, . 
	\label{condstret}
    \end{equation}
\end{corollary}

\begin{proof} \textit{Necessity.}
    We remark that if $E-k\D$ is $L^{p}$-dissipative, then
    \begin{equation}
	\cases{k\leq 1 & if $p=2$\cr
	k<1 & if $p\neq 2$.
	}
	\label{condp}
    \end{equation}
    In fact, in view of Theorem \ref{lemma:5}, we 
    have the necessary condition
    \beginrighe{condalgelas2}
    -(1-2/p)^{2}[(1-k)|\xi|^{2}+(1-2\n)^{-1}(\xi_{j}\o_{j})^{2}](\l_{j}\o_{j})^{2}
    \cr
    +
    (1-k)|\xi|^{2}|\l|^{2} + (1-2\n)^{-1}(\xi_{j}\l_{j})^{2}\geq 0
	\endrighe
    for any $\xi,\, \l,\, \o\in\R^{2}$, $|\o|=1$. If we take
    $\xi=(1,0)$, $\l=\o=(0,1)$ in \reff{condalgelas2} we find
    $$
    {4\over p\, p'}\, (1-k) \geq 0
    $$
    and then $k\leq 1$ for any $p$. If $p\neq 2$ and $k=1$, taking 
    $\xi=(1,0)$, $\l=(0,1)$, $\o=(1/\sqrt{2},1/\sqrt{2})$ in 
    \reff{condalgelas2},
    we find $-(1-2/p)^{2}(1-2\n)^{-1}\geq 0$. On the other hand,
    taking  $\xi=\l=(1,0)$, $\o=(0,1)$ we find $(1-2\n)^{-1}\geq 0$. 
    This is a contradiction and \reff{condp} is proved. 
    
   It is clear that if $E-k\D$ is $L^{p}$-dissipative, then
   $E-k'\D$ is $L^{p}$-dissipative for any $k'<k$. Therefore
   it is not restrictive to suppose that
   $E-k\D$ is $L^{p}$-dissipative for
   some $0<k<1$. Moreover $E$ is also $L^{p}$-dissipative.
   
    The $L^{p}$-dissipativity of $E-k\D$ ($0<k<1$) is 
    equivalent to the
    $L^{p}$-dissipativity of the operator 
   \begin{equation}
	E'u = \D u + (1-k)^{-1}(1-2\n)^{-1}\nabla \dive u.
       \label{defE'}
   \end{equation}
    
    Setting 
 \begin{equation}
	\n'=\n(1-k) +k/2,
     \label{defn'}
 \end{equation}
    we have
    $(1-k)(1-2\n)=1-2\n'$. 
    Theorem \ref{th:cnesdiss}  shows that 
  \begin{equation}
      {4\over p\, p'}\geq {1\over (3-4\n')^{2}}.
      \label{condelast'}
      \end{equation}
      
      Since $3-4\n'=3-4\n -2k(1-2\n)$, condition \reff{condelast'} 
      means
      $|3-4\n -2k(1-2\n)|\geq \sqrt{p\,p'}/2$, i.e.
      \begin{equation}
	  \left| k -{3-4\n\over 2(1-2\n)}\right| \geq
	  {\sqrt{p\,p'}\over 4|1-2\n|}
	  \label{condlastbis}
      \end{equation}

      Note that  the $L^{p}$-dissipativity of $E$ implies that
      \reff{condvecchia} holds. In particular we have
      $(3-4\n)/(1-2\n)> 0$. Hence \reff{condlastbis}
      is satisfied if either
      \begin{equation}
	  k \leq {1\over 2|1-2\n|}\left(|3-4\n| -
	  {\sqrt{p\,p'}\over 2}\right)
	  \label{bis1}
      \end{equation}
      or
      \begin{equation}
	  k \geq {1\over 2|1-2\n|}\left(|3-4\n| +
	  {\sqrt{p\,p'}\over 2}\right)
	  \label{bis2}
      \end{equation}
      
      Since  $$
      {|3-4\n|\over 2|1-2\n|}-1={3-4\n\over 2(1-2\n)}-1=
      {1\over 2(1-2\n)}\geq  
      -{\sqrt{p\,p'}\over 4|1-2\n|}
      $$
      we have
      $$
      {1\over 2|1-2\n|}\left(|3-4\n| +
		{\sqrt{p\,p'}\over 2}\right)\geq 1
      $$
      and \reff{bis2} is impossible. Then \reff{bis1} holds.
    Since $k>0$, we have the strict inequality in 
     \reff{condelast} and \reff{condstret} is proved.

      \textit{Sufficiency.} Suppose 
      \reff{condstret}. Since
      $$
      {4\over p\, p'}> {1\over (3-4\n)^{2}}\, ,
      $$
      we can take $k$ such that
     \begin{equation}
	  0< k <  {1\over 2|1-2\n|}\left(|3-4\n| -
	       {\sqrt{p\,p'}\over 2}\right).
	 \label{ksuchthat}
     \end{equation}
     
      Note that
      $$
	   {|3-4\n|\over 2|1-2\n|}-1={3-4\n\over 2(1-2\n)}-1=
	   {1\over 2(1-2\n)}\leq  
	   {\sqrt{p\,p'}\over 4|1-2\n|}.
	   $$
	   
	  This means
	   $$
	   {1\over 2|1-2\n|}\left(|3-4\n| -
		     {\sqrt{p\,p'}\over 2}\right)\leq 1
	   $$
	   and then $k<1$. Let $\n'$ be given by \reff{defn'}.
	   The $L^{p}$-dissipativity of $E-k\D$ is equivalent to the 
	   $L^{p}$-dissipativity of the operator $E'$ defined by
	   \reff{defE'}.

	   Condition \reff{condlastbis} (i.e. \reff{condelast'})
	   follows from 
	   \reff{ksuchthat} and Theorem \ref{th:cnesdiss} gives the result.
\end{proof}

\begin{corollary}
    There exists $k<2$ such that $k\D-E$ is $L^{p}$-dissipative 
    if and only if 
    \begin{equation}
	\left({1\over 2}-{1\over p}\right)^{2} < 
	      {2\n(2\n-1)\over
		    (1-4\n)^{2}}
	\, . 
	\label{eq:2nu}
    \end{equation}
\end{corollary}

\begin{proof}
    We may write $k\D-E=\widetilde{E}-\widetilde{k}\D$, where
    $\widetilde{k}=2-k$, 
    $\widetilde{E}=\D +(1-2\widetilde{\n})^{-1}\nabla\dive$, 
    $\widetilde{\n}=1-\n$.
    Theorem \ref{th:7} shows that $\widetilde{E}-\widetilde{k}\D$
    is $L^{p}$-dissipative if and only if
    \begin{equation}
	\left({1\over 2}-{1\over p}\right)^{2} < 
	     {2(\widetilde{\n}-1)(2\widetilde{\n}-1)\over
		   (3-4\widetilde{\n})^{2}}\, .
	    \label{eq:2nutilde}
	\end{equation}
	
	Condition \reff{eq:2nutilde} 
	coincides with\reff{eq:2nu} and the the Corollary is proved.
\end{proof}

\paragrafo{Dissipativity for a class of Systems of
Partial Differential Equations} \label{sec:PDE}

In this Section we consider a particular class of
operators \reff{eq:A}, namely the operators
\begin{equation}
    Au = \de_{h}(\A^{h}(x)\de_{h}u)
    \label{defspde}
\end{equation}
where $\A^{h}(x)=\{a_{ij}^{h}(x)\}$  ($i,j=1,\ldots,m$)
are matrices
with complex locally integrable entries defined in a domain 
$\O\subset\R^{n}$ ($h=1,\ldots,n$).

Our goal is to prove that $A$ is $L^{p}$-dissipative 
if and only if the algebraic condition
$$\display{
 \Re \lan \A^{h}(x) \l,\l\ran-(1-2/p)^{2}\Re\lan 
  \A^{h}(x)\o,\o\ran (\Re \lan\l,\o\ran)^{2}
   \cr
    -
    (1-2/p)\Re(\lan \A^{h}(x)\o,\l\ran -
    \lan \A^{h}(x)\l,\o\ran)
       \Re \lan \l,\o\ran   \geq 0
}
$$
   is satisfied for almost every $x\in\O$ and for
   every $\l,\o\in\C^{m}$, $|\o|=1$, $h=1,\ldots,n$. 
In order to obtain such a result, in
the next subsections we study the
dissipativity for some systems of ordinary differential
equations.

\subsection{Dissipativity for Systems of Ordinary Differential
Equations}

In this subsection we consider the operator $A$ defined by
\begin{equation}
    A u = (\A(x)u')'
    \label{defAord}
\end{equation}
where $\A(x)=\{a_{ij}(x)\}$ ($i,j=1,\ldots,m$) is a matrix with 
complex locally integrable entries defined in the bounded or 
unbounded
interval $(a,b)$.

In this case the sesquilinear form  $\elle(u,v)$ is given by
$$
\elle(u,v)=\int_{a}^{b}\lan\A u',v'\ran \, dx .
$$

\begin{lemma}\label{th:cond1ode}
    The operator $A$ is $L^{p}$-dissipative if and only if
    \beginrighe{cond1ode}
	  \int_{a}^{b}\Big(\Re \lan \A v',v'\ran-(1-2/p)^{2}|v|^{-4}\Re\lan \A v,v\ran 
	  (\Re \lan v,v'\ran)^{2} \cr
	  -(1-2/p)|v|^{-2}\Re(\lan \A v,v'\ran 
	  -\lan \A v',v\ran) \Re \lan v,v'\ran
	  \Big)\, dx  \geq 0
	 \endrighe
      for any $v\in (\Cspt^{1}((a,b)))^{m}$. 
\end{lemma}

\begin{proof}
    It is a particular case of Lemma \ref{lemma:5}.
\end{proof}

\begin{theorem}\label{th:0}
    The operator $A$ is $L^{p}$-dissipative if and only if
    \beginrighe{condode}
        \Re \lan \A(x) \l,\l\ran-(1-2/p)^{2}\Re\lan \A(x)\o,\o\ran (\Re \lan\l,\o\ran)^{2}
	\cr
	-
	(1-2/p)\Re(\lan \A(x)\o,\l\ran -\lan \A(x)\l,\o\ran)
	\Re \lan \l,\o\ran   \geq 0
    \endrighe
    for almost every $x\in(a,b)$ and for any 
    $\l,\o\in\C^{m}$, $|\o|=1$.
\end{theorem}

\begin{proof}  
    \textit{Necessity.}
   First we prove the result assuming that
   the coefficients $a_{ij}$ are constant and that $(a,b)=\R$.

    Let us fix $\l$ and $\o$ in $\C^{m}$, with $|\o|=1$, and choose
    $v(x)=\h(x/R)\, w(x)$ where
    $$
    w_{j}(x)=\cases{\m\o_{j} & if $x<0$\cr
    \m\o_{j}+x^{2}(3-2x)\l_{j}   & if $0\leq x\leq 1$\cr
    \m\o_{j}+\l_{j}    & if $x>1$,}
    $$
    $\m,R\in\R^{+}$, $\h\in\Cspt^{\infty}(\R)$, $\spt\h\subset 
    [-1,1]$ and $\h(x)=1$ if $|x|\leq 1/2$.

  We have 
   $$\display{
    \lan\A v',v'\ran  =\cr
    \lan\A w',w'\ran (\h(x/R))^{2}
   + 
   R^{-1}  (\lan\A w',w\ran +
  \lan \A w,w'\ran) \h(x/R)\h'(x/R) 
  \cr
  +R^{-2} \lan\A w,w\ran (\h'(x/R))^{2}
  }
   $$
    and then 
    $$
    \int_{\R}\lan\A v',v'\ran \, dx=
    \int_{0}^{1}\lan \A w',w'\ran\, dx
    + {1\over R^{2}}\int_{-R}^{R}  \lan\A w,w\ran (\h'(x/R))^{2}dx
    $$
    provided that $R>2$. Since $\lan\A w,w\ran$ is bounded, we have
    $$
    \lim_{R\to+\infty}\int_{\R}\lan\A v',v'\ran \, dx=
    \int_{0}^{1}\lan \A w',w'\ran\, dx.
    $$
    
    On the set where $v\neq 0$ we have 
    $$\display{|v|^{-4}\lan\A v,v\ran (\Re \lan 
    v,v'\ran)^{2} = |w|^{-4}\lan \A w,w\ran (\Re \lan w,w'\ran)^{2}
    (\h(x/R))^{2}\cr
    + 2\, R^{-1} |w|^{-2}\lan \A w,w\ran \Re \lan w,w'\ran
    \h(x/R)) \, \h'(x/R) +
     R^{-2}\lan \A w,w\ran ( \h'(x/R))^{2}
    }
    $$
    form which it follows
    $$
    \display{
	\lim_{R\to+\infty}\int_{\R}|v|^{-4}\lan\A v,v\ran (\Re \lan 
	v,v'\ran)^{2}dx =
	\int_{0}^{1}|w|^{-4}\lan\A w,w\ran (\Re \lan 
	    w,w'\ran)^{2}dx.
	    }
    $$
    
    In the same way we obtain
    $$\display{
	\lim_{R\to+\infty}\int_{\R}|v|^{-2}(\lan\A v,v'\ran -
	\lan\A v',v\ran) \Re \lan 
	   v,v'\ran \, dx =\cr
	   \int_{0}^{1}|w|^{-2}(\lan\A w,w'\ran -
	\lan\A w',w\ran) \Re \lan 
	   w,w'\ran \, dx.
     }	
    $$
    
    Since $v\in (\Cspt^{1}(\R))^{m}$, we can put $v$ in \reff{cond1ode}. 
    Letting $R\to+\infty$, we find
    \beginrighe{cond1odew}
	     \int_{0}^{1}
	     \Big(\Re \lan \A w',w'\ran-(1-2/p)^{2}|w|^{-4}\Re\lan \A w,w\ran 
	     (\Re \lan w,w'\ran)^{2} \cr
	     -(1-2/p)|w|^{-2}\Re(\lan \A w,w'\ran 
	     -\lan \A w',w\ran) \Re \lan w,w'\ran
	     \Big)\, dx  \geq 0.
	     \endrighe
    
	    On the interval $(0,1)$ we have
	     $$\display{
	     \lan \A w',w'\ran = \lan\A \l,\l\ran\, 36 x^{2}(1-x)^{2},\cr
	     |w|^{-4}\lan \A w,w\ran 
			  (\Re \lan w,w'\ran)^{2}=
			  |\m\o+x^{2}(3-2x)\l|^{-4}\times\cr
       (\m^{2}\lan\A\o,\o\ran 
       +\m(\lan\A\o,\l\ran+\lan\A\l,\o\ran)x^{2}(3-2x) 
       +\lan\A\l,\l\ran x^{4}(3-2x)^{2})\times \cr
      [\Re (\m\lan\o,\l\ran\, 6x(1-x) + |\l|^{2}6x^{3}(3-2x)(1-x))]^{2}, 
       \cr
       |w|^{-2}(\lan \A w,w'\ran 
		    -\lan \A w',w\ran) \Re \lan w,w'\ran =
		    |\m\o+x^{2}(3-2x)\l|^{-2}\times\cr
		   \m(\lan\A\o,\l\ran - \lan\A\l,\o\ran)\, 6x(1-x) 
		   \Re (\m\lan\o,\l\ran\, 6x(1-x) + |\l|^{2}6x^{3}(3-2x)(1-x)).
		    }
	     $$

	  Letting $\m\to\infty$ in \reff{cond1odew} we find
	   $$\display{
		     36 \int_{0}^{1}\Big(\Re \lan \A \l,\l\ran-(1-2/p)^{2}\Re\lan \A \o,\o\ran 
		      (\Re \lan \o,\l\ran)^{2} \cr
		      -(1-2/p)\Re(\lan \A \o,\l\ran 
		      -\lan \A \l,\o\ran) \Re \lan \o,\l\ran
		     \Big)\, x^{2}(1-x)^{2} dx  \geq 0
		     }
		     $$
    and \reff{condode} is proved. 
    
    If $a_{hk}$ are not necessarily constant, 
    consider 
    $$
    v(x)=\e^{-1/2}\p((x-x_{0})/\e)
    $$
    where $x_{0}$ is a fixed point in $(a,b)$,
    $\p\in(\Cspt^{1}(\R))^{m}$ and $\e$ is 
    sufficiently small.

    In this case \reff{cond1ode} shows that
       $$\display{
	   \int_{\R}\Big(\Re \lan \A (x_0+\e y)\p',\p'\ran 
	   -(1-2/p)^{2}|\p|^{-4}\Re\lan \A (x_0+\e y)\p,\p\ran (\Re \lan 
	   \p,\p'\ran)^{2} \cr
	   -(1-2/p)|\p|^{-2}\Re(\lan \A (x_0+\e y)\p,\p'\ran
   -\lan \A (x_0+\e y)\p',\p\ran) \Re \lan \p,\p'\ran
	   \Big)\, dy  \geq 0.
   }      
	   $$
    
	   Letting $\e\to 0^{+}$ we find for almost every $x_{0}$
	   $$\display{
		   \int_{\R}\Big(\Re \lan \A (x_0)\p',\p'\ran 
		   -(1-2/p)^{2}|\p|^{-4}\Re\lan \A (x_0)\p,\p\ran (\Re \lan 
		   \p,\p'\ran)^{2} \cr
		   -(1-2/p)|\p|^{-2}\Re(\lan \A (x_0)\p,\p'\ran
	   -\lan \A (x_0)\p',\p\ran) \Re \lan \p,\p'\ran
		   \Big)\, dy  \geq 0.
	   }      
		   $$
   
		  Because this inequality holds for any 
		  $\p\in\Cspt^{1}(\R)$, 
		 what we have obtained 
   for constant 
    coefficients gives the result.

    \textit{Sufficiency.} It is clear that, if \reff{condode} holds, 
    then the integrand in \reff{cond1ode} is nonnegative almost 
    everywhere and  Lemma \ref{th:cond1ode} gives the result.
\end{proof}

\begin{corollary}\label{cor:pos}
    If the operator $A$ is $L^{p}$-dissipative, then
    $$
    \Re \lan \A(x)\l,\l\ran \geq 0
    $$
    for almost every $x\in(a,b)$ and for any $\l\in\C^{m}$.
\end{corollary}

\begin{proof}
    Fix $x\in(a,b)$ such that \reff{condode} holds for any $\l,\o\in\C^{m}$,
    $|\o|=1$. For any $\l\in\C^{m}$, choose $\o$ such that $\lan\l, 
    \o\ran=0$, $|\o|=1$. The result follows by putting $\o$ in 
    \reff{condode}.
\end{proof}
¥

It is interesting to compare the operator $A$ with the
operator $I(d^{2}/dx^{2})$. 

\begin{corollary}\label{cor:1}
    There exists $k>0$ such that $A-kI(d^{2}/dx^{2})$ is 
       $L^{p}$-dissipative if and only if 
       \begin{equation}
	   \essinf_{(x,\l,\o)\in (a,b)\times\C^{m}\times\C^{m}
	   \atop |\l|=|\o|=1} P(x,\l,\o) >0
	   \label{Icondsuk}
       \end{equation}
       where 
       $$\display{
       P(x,\l,\o)=  \Re \lan \A(x) \l,\l\ran-(1-2/p)^{2}\Re\lan \A(x)\o,\o\ran (\Re \lan\l,\o\ran)^{2}
	   \cr
	   -
	   (1-2/p)\Re(\lan \A(x)\o,\l\ran -\lan \A(x)\l,\o\ran)
	   \Re \lan \l,\o\ran  .
	   }
	   $$
       There exists $k>0$ such that $kI(d^{2}/dx^{2})-A$ is 
	   $L^{p}$-dissipative if and only if 
	   \begin{equation}
	       \esssup_{(x,\l,\o)\in (a,b)\times\C^{m}\times\C^{m}
	   \atop |\l|=|\o|=1} P(x,\l,\o) <\infty .
	       \label{IIcondsuk}
	   \end{equation}
   \end{corollary}

   \begin{proof}
       In view of Theorem \ref{th:0}, $A-kI(d^{2}/dx^{2})$ is $L^{p}$-dissipative if 
       and only if
       $$
       P(x,\l,\o)-k|\l|^{2}+k(1-2/p)^{2}(\Re\lan \l,\o\ran)^{2}\geq 0
       $$
       for almost every $x\in(a,b)$ and for any $\l,\o\in\C^{m}$, 
       $|\o|=1$. Since
     \begin{equation}
	   |\l|^{2}-(1-2/p)^{2}(\Re\lan \l,\o\ran)^{2}\geq {4\over p\,p'}\, 
	     |\l|^{2},
	 \label{dislambda}
     \end{equation}
       we can find a positive $k$ such that this is true if 
       and only if
       \begin{equation}
	   \essinf_{(x,\l,\o)\in (a,b)\times\C^{m}\times\C^{m}
		   \atop |\o|=1}  {P(x,\l,\o)\over |\l|^{2}-(1-2/p)^{2}(\Re \lan 
	     \l,\o\ran)^{2}} >0. 
	     \label{preIcondsuk}
	 \end{equation} 
       
	 On the other hand, inequality \reff{dislambda} shows that
	\begin{equation}
	     {P(x,\l,\o)\over |\l|^{2}} \leq
		  {P(x,\l,\o)\over |\l|^{2}-(1-2/p)^{2}(\Re \lan 
		      \l,\o\ran)^{2}} \leq
		      {p\, p'\over 4}\, {P(x,\l,\o)\over |\l|^{2}}
	    \label{eq:3.5bis}
	\end{equation}
	 and then \reff{preIcondsuk} and \reff{Icondsuk} are equivalent.
	 
       In the same way the operator $kI(d^{2}/dx^{2})-A$ is 
	   $L^{p}$-dissipative if and only if 
	   $$
	     -P(x,\l,\o)+k|\l|^{2}-k(1-2/p)^{2}(\Re\lan \l,\o\ran)^{2}\geq 0
	   $$
	   for almost every $x\in(a,b)$ and for any $\l,\o\in\C^{m}$, 
	       $|\o|=1$. We can find a positive $k$ such that 
	       this is true if 
       and only if 
       $$
       \esssup_{(x,\l,\o)\in (a,b)\times\C^{m}\times\C^{m}
	       \atop |\o|=1}  {P(x,\l,\o)\over |\l|^{2}-(1-2/p)^{2}(\Re \lan 
		\l,\o\ran)^{2}} < \infty.
		$$
       
       This inequality is equivalent to \reff{IIcondsuk} because of
       \reff{eq:3.5bis}.
   \end{proof}

\begin{corollary} There exists $k\in\R$ such that $A-kI(d^{2}/dx^{2})$ is 
    $L^{p}$-dissipative if and only if 
$$
        \essinf_{(x,\l,\o)\in (a,b)\times\C^{m}\times\C^{m}
	\atop |\l|=|\o|=1} P(x,\l,\o) > -\infty.
$$
\end{corollary}

\begin{proof}
    The result can be proved as in Corollary \ref{cor:1}.
\end{proof}

\subsection{Real coefficient operators}

In the following we need the lemma
\begin{lemma}
    Let $0<\m_{1}\leq \m_{2}\leq \ldots \leq \m_{m}$. We have
   \begin{equation}
        \max_{\o\in\R^{m}\atop |\o|=1}\ [(\m_{h}\o_{h}^{2})(\m_{k}^{-1}\o_{k}^{2})] =
	{(\m_{1}+\m_{m})^{2}\over 4\, \m_{1}\m_{m}}.
       \label{exprmax}
   \end{equation}
\end{lemma}

\begin{proof}
    First we proof by induction on $m$ that
    \begin{equation}
	\max_{\o\in\R^{m}\atop |\o|=1}\ [(\m_{h}\o_{h}^{2})(\m_{k}^{-1}\o_{k}^{2})] =
	\max_{1\leq i < j\leq m}	{(\m_{i}+\m_{j})^{2}\over 
	4\, \m_{i}\m_{j}}.
        \label{exprmax2}
    \end{equation}
    In the case $m=2$, \reff{exprmax2} is equivalent to
    $$
    \max_{\f\in[0,2\pi]} [ \cos^{4}\f +\sin^{4}\f +
     (\m_{1}\m_{2}^{-1}+ \m_{2}\m_{1}^{-1})\cos^{2}\f\, \sin^{2}\f]
     = 	{(\m_{1}+\m_{2})^{2}\over 4\, \m_{1}\m_{2}},
    $$
    which can be easily proved.
    
   Let $m>2$ and suppose $\m_{1}<\m_{2}<\ldots<\m_{m}$; 
   the maximum of the left hand side of \reff{exprmax2}
   is the  maximum of the function
   $$
   \m_{h}\m_{k}^{-1}x_{h}x_{k}
   $$
   subject to the constraint $x\in K$, where $K=\{ x\in\R^{m}\ |\
   x_{1}+\ldots +x_{m}=1,\ 0\leq x_{j}
   \leq 1\ (j=1,\ldots,m)\}$. 
   To find the constrained maximum, we first examine the system
   \begin{equation}
       \cases{\g_{hk}x_{k}-\l=0  & $h=1,\ldots,m$\cr
       x_{1}+\ldots +x_{m}=1}
       \label{lagrange}
   \end{equation}
   with $0\leq x_{j}
   \leq 1\ (j=1,\ldots,m)$, where $\l$ is the Lagrange multiplier and $\g_{hk}=\m_{h}\m_{k}^{-1}+\m_{k}\m_{h}^{-1}$.
   
   Consider the homogeneous system
   \begin{equation}
       \g_{hk}x_{k}=0 \quad (h=1,\ldots,m).
       \label{omogeneo}
   \end{equation}
   
   One checks directly 
   that the vectors $x^{(k)}=(x_{1}^{(k)},\ldots,x_{m}^{(k)})$, 
   $$
   x_{1}^{(k)}={\m_{1}\over \m_{k}}\, {\m_{k}^{2}-\m_{2}^{2}\over
   \m_{2}^{2}-\m_{1}^{2}}, \quad
   x_{2}^{(k)}={\m_{2}\over \m_{k}}\, {\m_{1}^{2}-\m_{k}^{2}\over
   \m_{2}^{2}-\m_{1}^{2}}, \quad x_{j}^{(k)}=\d_{jk}\ (j=3,\ldots,m)
   $$
   for $k=3,\ldots,m$, are $m-2$ linearly independent eigensolutions of the system 
   \reff{omogeneo}. On the other hand,
   the determinant
   $$
   \left|\matrix{\g_{11} & \g_{12}\cr \g_{12} & \g_{22}}\right|=
   4-\g_{12}^{2}=-{(\m_{1}^{2}-\m_{2}^{2})^{2}\over 
   \m_{1}^{2}\m_{2}^{2}}<0
   $$
   and then the rank of the matrix $\{\g_{hk}\}$ is $2$.
   
   Therefore there exists a solution of the system
\begin{equation}
       \g_{hk}x_{k}= \l \qquad (h=1,\ldots,m)
    \label{sistemlam}
\end{equation}
   if and only if the vector $(\l,\ldots,\l)$ is orthogonal to
   any eigensolution  of 
   the adjoint homogeneous system. Since the matrix $\{\g_{hk}\}$ 
   is symmetric, there exists a solution of the system 
   \reff{sistemlam} if and only if
   \begin{equation}
       \l(x_{1}^{(k)}+ \cdots +x_{m}^{(k)})=0
       \label{condcomp}
   \end{equation}
   for $k=3,\ldots,m$.
   
   But
   $$
   x_{1}^{(k)}+ \cdots +x_{m}^{(k)}= - {\m_{1}\m_{2}+\m_{k}^{2}\over
   \m_{k}(\m_{1}+\m_{2})}+1=- {(\m_{k}-\m_{1})(\m_{k}-\m_{2})\over
   \m_{k}(\m_{1}+\m_{2})}<0
   $$
   and \reff{condcomp} are satisfied if and only if $\l=0$.
   This means that the system \reff{sistemlam} is solvable only 
   when $\l=0$ and the solutions are given by
   $$
   x =\sum_{k=3}^{m}u_{k}x^{(k)}
   $$
   for arbitrary $u_{k}\in\R$. On the other hand we are looking for 
   solutions of \reff{lagrange} with $0\leq x_{j}\leq 1$.
   Since $x_{j}=u_{j}$ for $j=3,\ldots,m$, we have $u_{j}\geq 0$. 
   This implies that
   $$
   x_{2}=\sum_{k=3}^{m}{\m_{2}\over \m_{k}}\, {\m_{1}^{2}-\m_{k}^{2}\over
   \m_{2}^{2}-\m_{1}^{2}}\, u_{k}\leq 0
   $$
   and since we require $x_{2}\geq 0$, we  have
   $u_{k}=0$ ($k=3,\ldots,m$), i.e. $x=0$. This solution does not 
   satisfy the last equation in \reff{lagrange}. This means that
   there
   are no extreme points belonging to the interior of $K$. The maximum is therefore  attained on the 
   boundary of $K$, where at least one of the $x_{j}$'s is zero.  
   This shows that if \reff{exprmax2} is true for $m-1$, then it is true
   also for $m$.
   
   We have proved \reff{exprmax2} assuming $0<\m_{1}<\ldots<\m_{m}$;
   in case $\m_{i}=\m_{j}$ for some $i,j$, it is obvious how to 
   obtain the result for $m$ from the one for $m-1$.
   
   Finally, let us show that
   \begin{equation}
       {(\m_{i}+\m_{j})^{2}\over 
	       4\, \m_{i}\m_{j}} \leq
	       {(\m_{1}+\m_{m})^{2}\over 4\, \m_{1}\m_{m}}
       \label{formulamax}
   \end{equation}
   for any $1\leq i,j \leq m$.   
   Set $\m_{j}=\a_{j}\m_{m}$ and suppose $i\leq j$. 
   We have $0< \a_{1}\leq \ldots 
   \leq \a_{m}=1$. Inequality \reff{formulamax} is equivalent to
   $$
   \a_{1}(\a_{i}+\a_{j})^{2}\leq \a_{i}\a_{j}(\a_{1}+1)^{2}
   $$
   i.e.
   $$
   \a_{1}\a_{i}(\a_{i}-\a_{j})+(\a_{1}\a_{j}-\a_{i})\a_{j}\leq 0
   $$
   and this is true, because $\a_{i}\leq \a_{j}$ and
   $\a_{1}\a_{j}\leq \a_{1}\leq \a_{i}$.
\end{proof}

\begin{theorem}\label{th:3}
  Let $\A$ be a real matrix $\{a_{hk}\}$ with $h,k=1,\ldots,m$.
 Let us suppose $\A=\A^{t}$ and  $\A\geq 0$
 (in the sense $\lan\A(x)\xi,\xi\ran\geq 0$, for almost every 
 $x\in(a,b)$ and for any $\xi\in\R^{m}$).
 The operator $A$ is 
    $L^{p}$-dissipative if and only if
    $$
    \left({1\over 2}-{1\over 
	p}\right)^{2} (\m_{1}(x)
	+\m_{m}(x))^{2} \leq \m_{1}(x)\m_{m}(x)
    $$
    almost everywhere, 
    where $\m_{1}(x)$ and $\m_{m}(x)$ are the smallest and the largest 
    eigenvalues of the matrix $\A(x)$ respectively. In the particular 
    case $m=2$ this condition is equivalent to
    $$
    \left({1\over 2}-{1\over 
	p}\right)^{2}(\tr \A(x))^{2}  \leq \det \A(x) 
    $$
    almost everywhere.
\end{theorem}

\begin{proof}
    From Theorem \ref{th:0} $A$ is $L^{p}$-dissipative 
    if and only if \reff{condode} holds for almost every $x\in(a,b)$ 
    and for any $\l,\o\in\C^{m}$, $|\o|=1$. We claim that 
    in the present case this condition is equivalent to
\begin{equation}
        \lan \A(x) \xi,\xi\ran-
        (1-2/p)^{2}\lan \A(x) \o,\o\ran (\lan \xi,\o\ran)^{2} \geq 0
    \label{condreale}
\end{equation}
    for almost every $x\in(a,b)$ and for 
    any $\xi, \o\in\R^{m}$, $|\o|=1$. 
    Indeed, it is obvious that if 
   $$
	    \lan \A(x) \l,\l\ran-
	    (1-2/p)^{2}\lan \A(x) \o,\o\ran (\Re\lan \l,\o\ran)^{2} 
	    \geq 0
$$
    for almost every $x\in(a,b)$ and for 
	any $\l, \o\in\C^{m}$, $|\o|=1$, then \reff{condreale} holds 
	for almost every $x\in(a,b)$ and for 
    any $\xi, \o\in\R^{m}$, $|\o|=1$. 
   Conversely, fix $x\in(a,b)$ and 
    suppose that \reff{condreale} holds for any 
    $\xi, \o\in\R^{m}$, $|\o|=1$.
    Let $Q$ be an orthogonal matrix such that $\A(x)=Q^{t}DQ$,
    $D$ being a diagonal matrix. 
    If we denote by 
    $\m_{j}$ the eigenvalues of $\A(x)$,
    we have
    $$\display{
    \lan \A(x) \l,\l\ran-
	    (1-2/p)^{2}\lan \A(x) \o,\o\ran (\Re\lan \l,\o\ran)^{2}
	    \cr
	   = \lan DQ \l,Q\l\ran-
	    (1-2/p)^{2}\lan DQ \o,Q\o\ran (\Re\lan Q\l,Q\o\ran)^{2}
	    \cr
	   =  \m_{j}|(Q\l)_{j}|^{2} -
	   (1-2/p)^{2}(\m_{j}|(Q\o)_{j}|^{2})(\Re\lan Q\l,Q\o\ran)^{2}
	    \cr
	   \geq \m_{j}|(Q\l)_{j}|^{2} -
		       (1-2/p)^{2}(\m_{j}|(Q\o)_{j}|^{2})
		       (|(Q\l)_{k}|\, |(Q\o)_{k}|)^{2}.
	    }
    $$
    The last expression is nonnegative because of 
    \reff{condreale} and the equivalence is proved. 
    
    Let us fix $x\in(a,b)$. We may write
       \reff{condreale} as
   \begin{equation}
	   (1-2/p)^{2}(\m_{h}\o_{h}^{2})(\xi_{k}\o_{k})^{2}
	   \leq \m_{j}\xi_{j}^{2}
       \label{condgamma}
   \end{equation}
   for any $\xi, \o\in\R^{m}$, $|\o|=1$. 
   Let us fix $\o\in\R^{m}$, $|\o|=1$;
   inequality \reff{condgamma} is true if and only if
   $$
   (1-2/p)^{2}(\m_{h}\o_{h}^{2})\sup_{\xi\in\R^{n}\atop \xi\neq 0}
   {(\xi_{k}\o_{k})^{2}\over
   \m_{j}\xi_{j}^{2}} \leq 1.
   $$
   
   We have
   $$
   \max_{\xi\in\R^{n}\atop \xi\neq 0}
      {(\xi_{k}\o_{k})^{2}\over
      \m_{j}\xi_{j}^{2}} = \m_{k}^{-1}\o_{k}^{2};
   $$
   in fact, by Cauchy's inequality, 
   we have $(\xi_{k}\o_{k})^{2}\leq (\m_{j}\xi_{j}^{2})
   (\m_{k}^{-1}\o_{k}^{2})$ for any $\xi\in\R^{m}$ and there 
   is equality if $\xi_{j}=\m_{j}^{-1}\o_{j}$.
   
   Therefore \reff{condgamma} is satisfied if and only if
   $$
   (1-2/p)^{2}(\m_{h}\o_{h}^{2})(\m_{k}^{-1}\o_{k}^{2})\leq 1
   $$
   for any $\o\in\R^{m}$, $|\o|=1$, and \reff{exprmax} shows that 
   this is true if and only if
   $$
   \left({1\over 2}-{1\over 
	p}\right)^{2} {(\m_{1}+\m_{m})^{2}\over \m_{1}\m_{m}}\leq 1\, .
   $$

The result for $m=2$ 
follows from the identities 
\begin{equation}
    \m_{1}(x)\m_{2}(x)=\det \A(x) , \quad \m_{1}(x)+\m_{2}(x)=\tr \A(x) .
    \label{identeigen}
\end{equation}
\end{proof}

\begin{corollary}\label{cor:2}
    Let $\A$ be a real and symmetric matrix. 
    Denote by $\m_{1}(x)$ and $\m_{m}(x)$ the 
       smallest and the largest eigenvalues of $\A(x)$ respectively.
       There exists $k>0$ such that 
       $A-kI(d^{2}/dx^{2})$ is 
       $L^{p}$-dissipative if and only if
     \begin{equation}
	       \essinf_{x\in(a,b)}
	     \left[(1+\sqrt{p\,p'}/2)\,\m_{1}(x)+ (1-\sqrt{p\,p'}/2)\,
	     \m_{m}(x)
	     \right]
		     >0.
	 \label{primacondmu}
     \end{equation}
   In the particular case  $m=2$ conditions \reff{primacondmu} is
   equivalent to
   \begin{equation} 
	   \essinf_{x\in(a,b)}\left[\tr \A(x) - 
	   {\sqrt{p\,p'}\over 2}\sqrt{(\tr\A(x))^{2}- 
	   4\det\A(x)} \right]
	   >0.
       \label{primacondmu2}
   \end{equation}
   \end{corollary}

   \begin{proof}
       \textit{Necessity.}
       Corollary \ref{cor:pos} shows that $\A(x)-kI\geq 0$ almost 
            everywhere.  
       In view of Theorem \ref{th:3}, we have that $A-kI(d^{2}/dx^{2})$ 
       is $L^{p}$-dissipative if and only if
   \begin{equation}
       \left({1\over p}-{1\over 
	      2}\right)^{2} (\m_{1}(x)
	      +\m_{m}(x)-2k)^{2} \leq (\m_{1}(x)-k)\, (\m_{m}(x)-k)
       \label{conddeterk}
   \end{equation}
       almost everywhere.

    Inequality \reff{conddeterk} is
    	  \begin{equation}
    	       {4\over p\,p'}\, (\m_{1}(x)
    	                  +\m_{m}(x)-2k)^{2} -  (\m_{1}(x)-\m_{m}(x))^{2}
    	                  \geq 0.
    	      \label{423}
    	  \end{equation}
    
    	By Corollary \ref{cor:1}, $A-k'I(d^{2}/dx^{2})$ is 
    	$L^{p}$-dissipative for any $k' \leq k$. Therefore 
    	inequality \reff{423} holds if we replace $k$ by any $k'< k$.
    	This implies that $k$ is less than or equal to the smallest root of
    	the left hand-side of \reff{423}, i.e.
    	\begin{equation}
    	       k \leq {1\over 2}\left[(1+\sqrt{p\,p'}/2)\,
    	       \m_{1}(x)+ (1-\sqrt{p\,p'}/2)\,
    			\m_{m}(x)
    			\right]
    		\label{primacondk}
    	    \end{equation}
	    and \reff{primacondmu} is proved.
	    
	\textit{Sufficiency.} 
	Let $k$ be such that
$$
        0<k \leq  \essinf_{x\in(a,b)} {1\over 2}
             \left[(1+\sqrt{p\,p'}/2)\,\m_{1}(x)+ (1-\sqrt{p\,p'}/2)\,
             \m_{m}(x)
             \right]
$$
    Since $\m_{1}(x)\leq 
     \m_{m}(x)$ and
     $\sqrt{p\,p'}/2\geq 1$, we have
  \begin{equation}
         (1+\sqrt{p\,p'}/2)\,\m_{1}(x)+ (1-\sqrt{p\,p'}/2)\,
                \m_{m}(x)\leq 
            2\, \m_{1}(x)
      \label{newprima}
  \end{equation}
     and then $\A(x)-kI\geq 0$ almost everywhere. 
     The constant $k$ satisifies \reff{primacondk}
    and this implies \reff{423}, i.e. \reff{conddeterk}.
     Theorem \ref{th:3} gives the result.
     
     The equivalence between \reff{primacondmu} and \reff{primacondmu2}
     follows from the identities \reff{identeigen}.
     \end{proof}
     
     If we require something more about the matrix $\A$ we have also
     \begin{corollary}
     		Let $\A$ be a real and symmetric matrix. 
     		Suppose $\A\geq 0$ almost everywhere.
     		Denote by $\m_{1}(x)$ and $\m_{m}(x)$ the 
     		   smallest and the largest eigenvalues of $\A(x)$ respectively.
     		 If  there exists $k>0$ such that 
     		   $A-kI(d^{2}/dx^{2})$ is 
     		   $L^{p}$-dissipative, then
     		 \begin{equation}
     			   \essinf_{x\in(a,b)}
     			 \left[
     			 \m_{1}(x)\m_{m}(x) - \left({1\over 2}-{1\over p}\right)^{2}
     			 (\m_{1}(x)+\m_{m}(x))^{2}
     			 \right]
     				 >0.
     		     \label{primacondmun}
     		 \end{equation}
		 If, in addition,  there exists
		$C$ such that 
		\begin{equation}
		    \lan\A(x)\xi,\xi\ran \leq C|\xi|^{2}
		    \label{condC}
		\end{equation}
		for almost every $x\in(a,b)$ and for any $\xi\in\R^{m}$,
		the converse is also true.
        In the particular case  $m=2$ condition \reff{primacondmun} 
        is equivalent to
        $$
             \essinf_{x\in(a,b)}\left[\det \A(x)
     	       - \left({1\over 2}-{1\over 
     	       p}\right)^{2}
     	    (\tr\A(x))^{2}\right] > 0 .
         $$
        \end{corollary}
     
     \begin{proof}
         \textit{Necessity.} By Corollary \ref{cor:2}, 
         \reff{primacondk} holds. On the other hand we have
	 	    $$\display{
	 	    \left[(1+\sqrt{p\,p'}/2)\,
	 			   \m_{1}(x)+ (1-\sqrt{p\,p'}/2)\,
	 				    \m_{m}(x)
	 				    \right]\cr
	 				    \leq
	         \left[(1-\sqrt{p\,p'}/2)\,\m_{1}(x)+ (1+\sqrt{p\,p'}/2)\,
	 					       \m_{m}(x) \right]
	 					       }
	 	    $$
	 	    and then
	 	    $$\display{
	 	    4k^{2}\leq 
	 	    \left[(1+\sqrt{p\,p'}/2)\,
	 				       \m_{1}(x)+ (1-\sqrt{p\,p'}/2)\,
	 						\m_{m}(x)
	 						\right]\cr
	 						\times
	 		    \left[(1-\sqrt{p\,p'}/2)\,\m_{1}(x)+ (1+\sqrt{p\,p'}/2)\,
	 			 \m_{m}(x) \right].
	 							   }
	 							   $$
	 	
	 	This inequality can be written as
$$
	         {4k^{2}\over p\,p'} \leq 
	              \m_{1}(x)\m_{m}(x) - \left({1\over 2}-{1\over p}\right)^{2}
	                           (\m_{1}(x)+\m_{m}(x))^{2}
$$
	 	and \reff{primacondmun} is proved.
		
		\textit{Sufficiency.} There exists $h>0$ such that
	 $$
	h \leq 
      \m_{1}(x)\m_{m}(x) - \left({1\over 2}-{1\over p}\right)^{2}
	(\m_{1}(x)+\m_{m}(x))^{2}
	 $$
	 almost everywhere, i.e.
	 $$\display{
	 p\, p' h \leq 
	 \left[(1+\sqrt{p\,p'}/2)\,
			\m_{1}(x)+ (1-\sqrt{p\,p'}/2)\,
			  \m_{m}(x) \right]\cr
			\times
		\left[(1-\sqrt{p\,p'}/2)\,\m_{1}(x)+ (1+\sqrt{p\,p'}/2)\,
			\m_{m}(x) \right]
			}
	 $$
	 almost everywhere. Since $\m_{1}(x)\geq 0$, we have also
	 \begin{equation}
	     (1-\sqrt{p\,p'}/2)\,\m_{1}(x)+ (1+\sqrt{p\,p'}/2)\,
	                      \m_{m}(x) \leq
	                      (1+\sqrt{p\,p'}/2)\,
	                   \m_{m}(x)    
	     \label{weha}
	 \end{equation}	 
	 and then
	 $$\display{
	 (1+\sqrt{p\,p'}/2)^{-1}p\, p' h \cr
	 \leq
	 \left[(1+\sqrt{p\,p'}/2)\,
				 \m_{1}(x)+ (1-\sqrt{p\,p'}/2)\,
				   \m_{m}(x) \right]
				   \esssup_{y\in(a,b)}\m_{m}(y)
				   }
	 $$
	 almost everywhere. By \reff{condC} $\esssup\m_{m}$ is finite and
	 by \reff{primacondmun} it is greater than zero. Then
	 \reff{primacondmu} holds and Corollary \ref{cor:2} gives the result.
     \end{proof}
     
    \begin{remark} Generally speaking, assumption \reff{condC} cannot be omitted,
	 even if $\A\geq 0$.	
	Consider, e.g., $(a,b)=(1,\infty)$,
	$m=2$,
	$\A(x)=\{a_{ij}(x)\}$ where
	$a_{11}(x)=(1-2/\sqrt{pp'})x + x^{-1}, \
	a_{12}(x)=a_{21}(x)=0, \
	a_{22}(x)=(1+2/\sqrt{pp'})x + x^{-1}$.
	We have
	$$\m_{1}(x)\m_{2}(x) - \left({1\over 2}-{1\over p}\right)^{2}
	(\m_{1}(x)+\m_{2}(x))^{2} = (8+4x^{-2})/(p\,p') 
	$$		
	and \reff{primacondmun} holds. But \reff{primacondmu} 
	is not satisfied, because
	$$
	(1+\sqrt{p\,p'}/2)\, \m_{1}(x)+ (1-\sqrt{p\,p'}/2)\,
					  \m_{2}(x) =2x^{-1} .
	$$
     \end{remark}
     
     \begin{corollary}
     Let $\A$ be a real and symmetric matrix. 
	 Denote by $\m_{1}(x)$ and $\m_{m}(x)$ the 
	    smallest and the largest eigenvalues of $\A(x)$ respectively.
	    There exists $k>0$ such that $kI(d^{2}/dx^{2})-A$ is 
		$L^{p}$-dissipative if and only if
		   \begin{equation} 
		       \esssup_{x\in(a,b)}
		    \left[(1-\sqrt{p\,p'}/2)\,\m_{1}(x)+ (1+\sqrt{p\,p'}/2)\,
		\m_{m}(x) \right]
		    <\infty.
	    \label{secondacondmu}
	\end{equation}
	In the particular case  $m=2$ condition
	\reff{secondacondmu}
	is equivalent to
	$$
		\esssup_{x\in(a,b)}\left[\tr \A(x) + 
		{\sqrt{p\,p'}\over 2}\sqrt{(\tr\A(x))^{2}- 
		    4\det\A(x)} \right]
		    <\infty.
	   $$
	\end{corollary}
     
     \begin{proof}
	 The proof runs as in Corollary \ref{cor:2}.
    We have that 
     $kI(d^{2}/dx^{2})-A$ is $L^{p}$-dissipative if and 
     only if \reff{conddeterk} holds, provided that
   $$
	 kI-\A(x)\geq 0
   $$
   almost everywhere.
     Because of this inequality, 
     we have to replace \reff{primacondk} and \reff{newprima} by
     $$
	     k \geq {1\over 2} 
	     \left[(1-\sqrt{p\,p'}/2)\,\m_{1}(x)+ (1+\sqrt{p\,p'}/2)\,
	       \m_{m}(x) \right]
      $$
	 and 
	 \begin{equation}
	 (1-\sqrt{p\,p'}/2)\,\m_{1}(x)+ (1+\sqrt{p\,p'}/2)\,
			   \m_{m}(x) 
			 \geq 2\, \m_{m}(x)
			 \label{seconcondk}
		     \end{equation}
	 respectively.
   \end{proof}

In the case of a positive matrix $\A$, we have

\begin{corollary}
   Let $\A$ be a real and symmetric matrix. Suppose $\A\geq 0$ almost everywhere.
       Denote by $\m_{1}(x)$ and $\m_{m}(x)$ the 
	  smallest and the largest eigenvalues of $\A(x)$ respectively.
	  There exists $k>0$ such that $kI(d^{2}/dx^{2})-A$ is 
	      $L^{p}$-dissipative if and only if
		 \begin{equation} 
		     \esssup_{x\in(a,b)}\m_{m}(x)
		  <\infty.
	  \label{supmu}
      \end{equation}
      \end{corollary}
      
      \begin{proof}
	  The equivalence between \reff{secondacondmu} and \reff{supmu}
	  follows 
          from \reff{weha} and \reff{seconcondk}.
      \end{proof}
      
We have also

\begin{corollary}  Let $\A$ be a real and symmetric matrix. 
    Denote by $\m_{1}(x)$ and $\m_{m}(x)$ the 
       smallest and the largest eigenvalues of $\A(x)$ respectively.
       There exists $k\in\R$ such that 
       $A-kI(d^{2}/dx^{2})$ is 
       $L^{p}$-dissipative if and only if
     $$
	       \essinf_{x\in(a,b)}
	     \left[(1+\sqrt{p\,p'}/2)\,\m_{1}(x)+ (1-\sqrt{p\,p'}/2)\,
	     \m_{m}(x)
	     \right]
		     > -\infty.
          $$
   In the particular case  $m=2$ this condition is
   equivalent to
   $$
	   \essinf_{x\in(a,b)}\left[\tr \A(x) - 
	   {\sqrt{p\,p'}\over 2}\sqrt{(\tr\A(x))^{2}- 
	   4\det\A(x)} \right]
	   >-\infty.
    $$
   \end{corollary}

   \begin{proof}
   The proof is similar to that of Corollary \ref{cor:2}. 
   \end{proof}

\subsection{$L^{p}$-dissipativity of the operator \reff{defspde}}

In this Section we consider the
partial differential operator \reff{defspde}
with complex coefficients.

Here $y_{h}$ denotes the $(n-1)$-dimensional
vector $(x_{1},\ldots,x_{h-1}, x_{h+1}, \ldots ,x_{n})$
and we set
$\o(y_{h})=\{ x_{h}\in\R\ |\ x\in\O\}$.

\begin{lemma}\label{lemma:5bis}
    The operator \reff{defspde} is $L^{p}$-dissipative if and only if
    the ordinary differential operators
    $$
        A(y_{h})[u(x_{h})]=d(\A^{h}(x)du/dx_{h})/dx_{h}
      $$
    are $L^{p}$-dissipative in $\o(y_{h})$
    for almost every $y_{h}\in \R^{n-1}$
    ($h=1,\ldots,n$). This condition is void if 
    $\o(y_{h})=\emptyset$.
\end{lemma}

\begin{proof} 
    \textit{Sufficiency.}
    Suppose $p\geq 2$. If $u\in (\Cspt^{1}(\O))^{m}$ we may write
    $$\display{
    \Re \sum_{h=1}^{n}\int_{\O}\lan \A^{h}(x)\de_{h}u,
    \de_{h}(|u|^{p-2}u)\ran dx=\cr
    \Re  \sum_{h=1}^{n}\int_{\R^{n-1}}dy_{h}\int_{\o(y_{h})}
    \lan \A^{h}(x)\de_{h}u, \de_{h}(|u|^{p-2}u)\ran dx_{h}.
    }
    $$
    By assumption
    $$
    \Re \int_{\o(y_{h})}
    \lan \A^{h}(x)v'(x_{h}), 
    (|v(x_{h})|^{p-2}v(x_{h}))'\ran dx_{h} \geq 0
    $$
   for almost every $y_{h}\in \R^{n-1}$ and
   for any $v\in(\Cspt^{1}(\o(y_{h})))^{m}$,
   provided $\o(y_{h})
   \neq \emptyset$ ($h=1,\ldots,n$). This implies
    $$
    \Re \sum_{h=1}^{n}\int_{\O}\lan \A^{h}(x)\de_{h}u,
       \de_{h}(|u|^{p-2}u)\ran dx \geq 0.
    $$
    
    The proof for $1<p<2$ runs in the same way. We have just to use
    \reff{eq:defdis2} instead of \reff{eq:defdis1}.
    
    \textit{Necessity.} Assume first that $\A^{h}$ are constant matrices
    and $\O=\R^{n}$. Let $p\geq 2$ and fix 
    $1\leq k \leq n$.
    
    Take $\a\in(\Cspt^{1}(\R))^{m}$ and
	$\b \in \Cspt^{1}(\R^{n-1})$. Consider
	$$
	u_{\e}(x)=\a(x_{k}/\e)\, \b(y_{k})
	$$

	We have
	 $$\display{
	 \sum_{h=1}^{n}\int_{\R^{n}}\lan \A^{h}\de_{h}u_{\e},\de_{h}
	 (|u_{\e}|^{p-2}u_{\e})
	 \ran dx = \cr
	 \e^{-2}\int_{\R^{n-1}}|\b(y_{k})|^{p}dy_{k}
	 \int_{\R}\lan \A^{k}
	 \a'(x_{k}/\e), 
	 \g'(x_{k}/\e)\ran 
	 \, dx_{k}
	 \cr
	 + \sum_{h=1\atop h\neq k}^{n}\int_{\R^{n-1}}\de_{h}\b(y_{k})\,
	 \de_{h}(|\b(y_{k})|^{p-2}\b(y_{k}))\, dy_{k}
	 \cr
	 \times
	 \int_{\R}\lan \A^{h}\a(x_{k}/\e),\a(x_{k}/\e)\ran \, |\a(x_{k}/\e)|^{p-2}dx_{k}
	 \cr
	 = 
	 \e^{-1}\int_{\R^{n-1}}|\b(y_{k})|^{p}dy_{k}
		  \int_{\R}\lan \A^{k}
		  \a'(t), 
		  (|\a(t)|^{p-2}\a(t))'\ran 
		  \, dt
		  \cr
		  + \e \sum_{h=1\atop h\neq 
		  k}^{n}\int_{\R^{n-1}}\de_{h}\b(y_{k})\,
		  \de_{h}(|\b(y_{k})|^{p-2}\b(y_{k}))\, dy_{k}
		  \int_{\R}\lan \A^{h}\a(t),\a(t)\ran 
		  \, |\a(t)|^{p-2}dt
	 }
      $$
      where $\g(t)=|\a(t)|^{p-2}\a(t)$.
    Keeping in mind 
    \reff{eq:defdis1} and letting $\e \to 0^{+}$, we find
    $$
    \Re \int_{\R^{n-1}}|\b(y_{k})|^{p}dy_{k}
		  \int_{\R}\lan \A^{k}
		  \a'(t), 
		  (|\a(t)|^{p-2}\a(t))'\ran 
		  \, dt\geq 0
    $$
    and then
    $$
    \Re  \int_{\R}\lan \A^{k}
		  \a'(t), 
		  (|\a(t)|^{p-2}\a(t))'\ran 
		  \, dt\geq 0
    $$
    for any $\a\in\Cspt^{1}(\R)$.
    This shows that $A(y_{k})$ is $L^{p}$-dissipative.
    
    If $\A^{h}$ are not necessarily constant, consider
    $$
    v(x)=\e^{(2-n)/2}\p((x-x_{0})/\e)
    $$
    where $x_{0}\in\O$, $\p\in(\Cspt^{1}(\R^{n}))^{m}$ and
    $\e$ is sufficiently small.

    In view of Lemma \ref{lemma:5} we write
    $$\display{
       \int_{\O}\Big(\Re \lan \A^{h} \de_{h}v,\de_{h}v\ran
       -(1-2/p)^{2}|v|^{-4}\Re\lan \A^{h} v,v\ran (\Re \lan 
       v,\de_{h}v\ran)^{2} \cr
       -(1-2/p)|v|^{-2}\Re(\lan \A^{h} v,\de_{h}v\ran
      - \lan \A^{h} \de_{h}v,v\ran)\Re \lan v,\de_{h}v\ran
	\Big) dx 
       \geq 0
       }
     $$
    i.e.
    $$\display{
	  \int_{\R^{n}}\Big(\Re \lan \A^{h}(x_{0}+\e z) \de_{h}\p,\de_{h}\p\ran
	  \cr
	  -(1-2/p)^{2}|\p|^{-4}\Re\lan \A^{h}(x_{0}+\e z) \p,\p\ran (\Re \lan 
	  \p,\de_{h}\p\ran)^{2} \cr
	  -(1-2/p)|\p|^{-2}\Re(\lan \A^{h}(x_{0}+\e z) \p,\de_{h}\p\ran
	  \cr
	 - \lan \A^{h}(x_{0}+\e z) \de_{h}\p,\p\ran)\Re \lan \p,\de_{h}\p\ran
	   \Big) dz 
	  \geq 0.
	  }
	$$
    
	Letting $\e\to 0^{+}$, we obtain
	$$\display{
		 \int_{\R^{n}}\Big(\Re \lan \A^{h}(x_{0}) \de_{h}\p,\de_{h}\p\ran
		 -(1-2/p)^{2}|\p|^{-4}\Re\lan \A^{h}(x_{0}) \p,\p\ran (\Re \lan 
		 \p,\de_{h}\p\ran)^{2} \cr
		 -(1-2/p)|\p|^{-2}\Re(\lan \A^{h}(x_{0}) \p,\de_{h}\p\ran
		- \lan \A^{h}(x_{0}) \de_{h}\p,\p\ran)\Re \lan \p,\de_{h}\p\ran
		  \Big) dy 
		 \geq 0
		 }
	       $$
	for almost every $x_{0}\in\O$. 
	
	Because of the arbitrariness of 
	$\p\in(\Cspt^{1}(\R^{n}))^{m}$, Lemma \ref{lemma:5} shows that the constant
	coefficient operator $\de_{h}(\A^{h}(x_{0})\de_{h})$ is 
	$L^{p}$-dissipative. From what has already been proved, 
	the ordinary differential operators $(\A^{h}(x_{0})v')'$ are
	$L^{p}$-dissipative ($h=1,\ldots,n$).

	Theorem \ref{th:0} yelds
	\beginrighe{condnew}
	 \Re \lan \A^{h}(x_{0}) \l,\l\ran-(1-2/p)^{2}\Re\lan 
	  \A^{h}(x_{0})\o,\o\ran (\Re \lan\l,\o\ran)^{2}
	   \cr
	    -
	    (1-2/p)\Re(\lan \A^{h}(x_{0})\o,\l\ran -
	    \lan \A^{h}(x_{0})\l,\o\ran)
	       \Re \lan \l,\o\ran   \geq 0
	   \endrighe
	for any $\l,\o\in\C^{m}$, $|\o|=1$, $h=1,\ldots,n$.
	
	Fix $h$ and denote by $N$ the set of $x_{0}\in\O$ such that
	\reff{condnew} does not hold for any 
	$\l,\o\in\C^{m}$, $|\o|=1$. Since $N$ has zero measure,
	for almost every $y_{h}\in \R^{n-1}$, 
	the cross-sections $\{x_{h}\in \R\ |\ 
	x\in N\}$ are measurable and have zero measure. 
	
	Hence, for almost every $y_{h}\in \R^{n-1}$, we
	have
	$$\display{
	 \Re \lan \A^{h}(x) \l,\l\ran-(1-2/p)^{2}\Re\lan 
	  \A^{h}(x)\o,\o\ran (\Re \lan\l,\o\ran)^{2}
	   \cr
	    -
	    (1-2/p)\Re(\lan \A^{h}(x)\o,\l\ran -
	    \lan \A^{h}(x)\l,\o\ran)
	       \Re \lan \l,\o\ran   \geq 0
	  }$$
	for almost every $x_{h}\in \o(y_{h})$ and 
	for any $\l,\o\in\C^{m}$, $|\o|=1$, provided
	$\o(y_{h}) \neq \emptyset$.
	The conclusion follows from Theorem \ref{th:0}.
    
    In the same manner we obtain the result for $1<p<2$.
\end{proof}

\begin{theorem}\label{thx:3}
    The operator \reff{defspde} is $L^{p}$-dissipative 
    if and only if
    \reff{condnew} holds for almost every $x_{0}\in\O$ and for any $\l,\o\in\C^{m}$, $|\o|=1$,
    $h=1,\ldots,n$.
\end{theorem}

\begin{proof}
    \textit{Necessity.} This has been already proved in the necessity
    part of the proof of Lemma \ref{lemma:5bis}. 
    
    \textit{Sufficiency.} We have seen that 
    if \reff{condnew} holds for
    almost every $x_{0}\in\O$ and
    for any $\l,\o\in\C^{m}$, $|\o|=1$,  
     the ordinary differential
    operator $A(y_{h})$ is $L^{p}$-dissipative 
    for almost every $y_{h}\in \R^{n-1}$, provided
    $\o(y_{h})\neq \emptyset$
    ($h=1,\ldots,n$). By Lemma \ref{lemma:5bis},  $A$ is 
    $L^{p}$-dissipative.    
\end{proof}

\begin{remark} In the scalar case ($m=1$), operator \reff{defspde} 
falls into the operators considered in  \cite{cialmaz}. 
In fact, if $Au=\sum_{h=1}^{n}\de_{h}(a^{h}\de_{h}u)$,
$a^{h}$ being a scalar function, $A$ can be written in the 
form \reff{cialmazop} with $\A=\{c_{hk}\}$,
$c_{hh}=a^{h}$, $c_{hk}=0$ if $h\neq k$.
The conditions obtained there can be directly compared with
\reff{condnew}.
The results of \cite{cialmaz} show that operator $A$
is $L^{p}$-dissipative if and only if \reff{cialmaz} holds.
This means that
\begin{equation}
    {4\over p\,p'}\,\lan \Re\A \xi,\xi\ran +
    \lan \Re\A \h,\h\ran -2(1-2/p)\lan \Im \A\xi,\h\ran \geq 0
    \label{vecchiacond}
\end{equation}
almost everywhere and for any $\xi,\h\in\R^{n}$ (see \cite[Remark 1, p.1082]
{cialmaz}). In this particular case
\reff{vecchiacond} is clearly equivalent to the following $n$ conditions
\begin{equation}
    {4\over p\,p'}\, (\Re a^{h})\, \xi^{2} + (\Re a^{h})\, \h^{2}
     -2(1-2/p)(\Im a^{h})\, \xi\h \geq 0
    \label{vecchiecond}
\end{equation}
almost everywhere and for any $\xi,\h\in \R$, $h=1,\ldots,n$.
On the other hand, in this case, \reff{condnew} reads as
\beginrighe{nuovecond}
    (\Re a^{h})|\l|^{2}-(1-2/p)^{2}(\Re a^{h})
    (\Re (\l\overline{\o})^{2}\cr
    -2(1-2/p)(\Im a^{h})\Re(\l\overline{\o})
   \Im(\l\overline{\o}) \geq 0
\endrighe
almost everywhere and for any $\l,\o\in\C$, $|\o|=1$, $h=1,\ldots,n$. Setting
$\xi+i\h=\l\overline{\o}$ and observing that
$|\l|^{2}=|\l\overline{\o}|^{2}=(\Re(\l\overline{\o}))^{2}+
(\Im(\l\overline{\o}))^{2}$, we see that conditions \reff{vecchiecond}
(and then \reff{vecchiacond}) are equivalent
to \reff{nuovecond}.
\end{remark}

In the case of a real coefficient operator \reff{defspde}, we have also

\begin{theorem}
    Let $A$ be the operator \reff{defspde}, where
$\A^{h}$ are real matrices $\{a^{h}_{ij}\}$ with $i,j=1,\ldots,m$.
 Let us suppose $\A^{h}=(\A^{h})^{t}$ and  $\A^{h}\geq 0$
 ($h=1,\ldots,n$).
 The operator $A$ is 
    $L^{p}$-dissipative if and only if
   \begin{equation}
	\left({1\over 2}-{1\over 
	   p}\right)^{2} (\m_{1}^{h}(x)
	   +\m_{m}^{h}(x))^{2} \leq 
	   \m_{1}^{h}(x)\,
	   \m_{m}^{h}(x)
       \label{ineigen}
   \end{equation}
    for almost every $x\in\O$, $h=1,\ldots,n$,  
    where $\m_{1}^{h}(x)$ and $\m_{m}^{h}(x)$ are the 
    smallest and the largest 
    eigenvalues of the matrix $\A^{h}(x)$ respectively. 
    In the particular 
    case $m=2$ this condition is equivalent to
    $$
    \left({1\over 2}-{1\over 
    p}\right)^{2}(\tr \A^{h}(x))^{2}  \leq \det 
    \A^{h}(x) 
    $$
    for almost every $x\in\O$, $h=1,\ldots,n$.
\end{theorem}

\begin{proof}
    By Theorem \ref{thx:3}, $A$ is $L^{p}$-dissipative if and only 
    if 
    $$
    \lan \A^{h}(x) \l,\l\ran-
    (1-2/p)^{2}\lan \A^{h}(x) \o,\o\ran (\Re\lan \l,\o\ran)^{2} \geq 0
 $$
  for almost every $x\in\O$, for any $\l,\o\in\C^{m}$, $|\o|=1$, 
  $h=1,\ldots,n$.  
    The proof of Theorem \ref{th:3} shows that
  these conditions are equivalent to \reff{ineigen}.
\end{proof}

\subsection{The angle of dissipativity}

In this Section we find the precise 
 angle of dissipativity for operator \reff{defspde}
 with complex coefficients.

We first consider the ordinary differential operator \reff{defAord} 
 where $\A(x)$ is a matrix whose elements are complex 
 locally integrable functions. 
 Define the functions
 \beginrighe{defPQ}
 P(x,\l,\o)=
 \Re \lan \A  \l,\l\ran-(1-2/p)^{2}\Re\lan \A \o,\o\ran (\Re \lan\l,\o\ran)^{2}
	\cr
	-
	(1-2/p)\Re(\lan \A \o,\l\ran -\lan \A \l,\o\ran)
	\Re \lan \l,\o\ran ; 
 \cr\cr
	Q(x,\l,\o)=\Im \lan \A  \l,\l\ran-(1-2/p)^{2}\Im\lan \A \o,\o\ran 
	(\Re \lan\l,\o\ran)^{2}
	\cr
	-
	(1-2/p)\Im(\lan \A \o,\l\ran -\lan \A \l,\o\ran)
	\Re \lan \l,\o\ran 
  \endrighe
 and denote by $\Xi$  the set 
 $$
 \Xi=\{ (x,\l,\o)\in (a,b)\times \C^{m}\times \C^{m}\ |\ 
 |\o|=1,\ P^{2}(x,\l,\o)+Q^{2}(x,\l,\o)> 0\}.
 $$

 By adopting the conventions introduced in Lemma \ref{lemma:elem},
 we have

 \begin{theorem}\label{th:4}
     Let $A$ be $L^{p}$-dissipative. The operator $zA$ is 
     $L^{p}$-dissipative if and only if
 $$
	 \th_{-}\leq \arg z \leq \th_{+}
 $$
     where 
     $$\display{
     \th_{-}=\arccot\left(\essinf_{(x,\l,\o)\in\Xi}
     (Q(x,\l,\o)/P(x,\l,\o))\right) - \pi, \cr
     \th_{+}=\arccot\left(\esssup_{(x,\l,\o)\in\Xi}
     (Q(x,\l,\o)/P(x,\l,\o))\right).
	 }
     $$
 \end{theorem}

 \begin{proof}
     In view of Theorem \ref{th:0} the operator $e^{i\th}A$ is 
     $L^{p}$-dissipative if and only if 
     \beginrighe{eith}
	    \Re \lan e^{i\th}\A  \l,\l\ran -(1-2/p)^{2}\Re\lan e^{i\th}\A \o,\o\ran (\Re \lan\l,\o\ran)^{2}
		 \cr
		 -
		 (1-2/p)\Re(\lan e^{i\th}\A \o,\l\ran -\lan e^{i\th}\A \l,\o\ran)
		 \Re \lan \l,\o\ran  
		 \geq 0
     \endrighe
     for almost every $x\in(a,b)$ and for any $\l,\o\in\C^{m}$, 
     $|\o|=1$.
     
     By means of the functions $P(x,\l,\o)$ and $Q(x,\l,\o)$ 
     introduced in \reff{defPQ}, we can write \reff{eith} in the form
 $$
	P(x,\l,\o)\, \cos \th   - Q(x,\l,\o)\, \sin \th  \geq 0.
 $$

 Lemma \ref{lemma:elem} gives the result.
 \end{proof}

 Let now $A$ be the partial differential operator \reff{defspde}. We have

\begin{theorem}
   Let $A$ be $L^{p}$-dissipative. The 
   operator $zA$ is $L^{p}$-dissipative if and only if
   $\th_{-}\leq \arg z \leq \th_{+}$, where
   $$\display{
   \th_{-}=\max_{h=1,\ldots,n}\
   \arccot\left(\essinf_{(x,\l,\o)\in\Xi_{h}}
       (Q_{h}(x,\l,\o)/P_{h}(x,\l,\o))\right) - \pi, 
   \cr
   \th_{+}=\min_{h=1,\ldots,n}\
   \arccot\left(\esssup_{(x,\l,\o)\in\Xi_{h}}
       (Q_{h}(x,\l,\o)/P_{h}(x,\l,\o))\right), 
   }
   $$
   and
   $$\display{
   P_{h}(x,\l,\o)=
\Re \lan \A^{h}(x) \l,\l\ran-(1-2/p)^{2}\Re\lan \A^{h}(x)
\o,\o\ran (\Re \lan\l,\o\ran)^{2}
       \cr
       -
       (1-2/p)\Re(\lan \A^{h}(x)\o,\l\ran -\lan \A^{h}(x)\l,\o\ran)
       \Re \lan \l,\o\ran , 
\cr\cr
       Q_{h}(x,\l,\o)=\Im \lan \A^{h}(x) \l,\l\ran-(1-2/p)^{2}
       \Im\lan \A^{h}(x)\o,\o\ran 
       (\Re \lan\l,\o\ran)^{2}
       \cr
       -
       (1-2/p)\Im(\lan \A^{h}(x)\o,\l\ran -\lan \A^{h}(x)\l,\o\ran)
       \Re \lan \l,\o\ran ,
   \cr\cr
   \Xi_{h}=\cr
   \{ (x,\l,\o)\in \O
   \times \C^{m}\times \C^{m}\ |\ 
|\o|=1,\ P_{h}^{2}(x,\l,\o)+Q_{h}^{2}(x,\l,\o)> 0\}.
   }
   $$
\end{theorem}

\begin{proof}
    By Theorem \ref{thx:3}, the operator
    $e^{i\th}A$ is $L^{p}$-dissipative if and only if
    \beginrighe{eith2}
	   \Re \lan e^{i\th}\A^{h}(x)  \l,\l\ran -(1-2/p)^{2}
	   \Re\lan e^{i\th}\A^{h}(x) \o,\o\ran (\Re \lan\l,\o\ran)^{2}
		\cr
		-
		(1-2/p)\Re(\lan e^{i\th}\A^{h}(x) 
		\o,\l\ran -\lan e^{i\th}\A^{h}(x) \l,\o\ran)
		\Re \lan \l,\o\ran  
		\geq 0
    \endrighe
    for almost every $x\in\O$ and for any $\l,\o\in\C^{m}$, 
    $|\o|=1$, $h=1,\ldots,n$.
   
    As in the proof of Theorems \ref{th:4}, conditions \reff{eith2}
    mean $\th_{-}^{(h)}\leq \th \leq \th_{+}^{(h)}$, where
    $$\display{
     \th_{-}^{(h)}=
     \arccot\left(\essinf_{(x,\l,\o)\in\Xi_{h}}
	 (Q_{h}(x,\l,\o)/P_{h}(x,\l,\o))\right) - \pi, 
     \cr
     \th_{+}^{(h)}=
     \arccot\left(\esssup_{(x,\l,\o)\in\Xi_{h}}
	 (Q_{h}(x,\l,\o)/P_{h}(x,\l,\o))\right), 
     }
     $$
     and the result follows.
\end{proof}

\end{document}